\documentclass[]{article}
\usepackage[misc]{ifsym}
\usepackage{amsmath}

\DeclareMathOperator*{\argmin}{arg\,min}
\usepackage{graphicx}
\usepackage{graphics}
\usepackage{amssymb,amsthm,epsfig}
\newtheorem{theorem}{Theorem}
\newtheorem{lemma}{Lemma}
\newtheorem{corollary}{Corollary}

\newtheorem{proposition}{Proposition}

\newtheorem{remark}{Remark}
\newcommand{\R}{\mathbb R}

\newcounter{reh}
\newcounter{rek}
\begin{document}
	\begin{center}
		{\large {\bf A Distribution Free   Truncated  Kernel Ridge Regression Estimator and Related Spectral Analyses}}\\
		\vskip 1cm  Asma Ben Saber$^a$ and  Abderrazek Karoui$^a$\footnote{
			\hspace*{-0.5cm} Corresponding Author: \hspace*{0.2cm} Abderrazek Karoui (abderrazek.karoui@fsb.ucar.tn)\\
			This work was supported in part by the  
			DGRST  research grant  LR21ES10 and the PHC-Utique research project 20G1503.
		} 
	\end{center}
	\vskip 0.5cm {\small
		
		\noindent $^a$ University of Carthage,
		Department of Mathematics, Faculty of Sciences of Bizerte, Tunisia.
	}\\

	\noindent{\bf Abstract}--- It is well known that kernel ridge regression (KRR)  is a   popular approach for  nonparametric regression based estimator.  Nonetheless, in the presence of a large data set with size $n\gg 1,$ the KRR based estimator has the  drawback to require an intensive  computational load.  Recently, some scalable KRR approaches  have been proposed with the aims to reduce the computational complexity of the KRR, while maintaining its superb convergence rate.  In this work, we study a new scalable KRR based
	approach for nonparametric regression. The proposed truncated kernel ridge regression (TKRR) approach is  simple. It  is based on substituting the full $n\times n$ random  kernel or Gram matrix $B_n,$ associated with a Mercer's kernel $\mathbb K,$ by its main $n\times N$ sub-matrix $A_N,$ where usually $N \ll n.$ Unlike some of the existing scalable KRR estimators, the TKRR  does not require a pre-processing step, such as an SVD or a sketching by a certain random matrix. Also,  we show that the TKRR works with $d-$dimensional random sampling  data  following an unknown probability law. To do so, we give a spectral analysis for  the compact kernel integral operator, associated with  a probability measure, different from its usual probability measure. This decay estimate is then extended to the decay of the  tail of the trace of the associated random Gram matrix. We expect that this spectral analysis is also useful for the study of others machine learning subjects such as the kernels based unsupervised domain  learning problems.   A special interest is devoted to develop rules for the optimal choices of the involved truncation order $N$ and  the value for regularization parameter $\lambda >0.$ The proposed  rules are based on the behavior and the decay rate of the spectrum of the positive integral operator, associated with the kernel $\mathbb K.$    These optimal values of the parameters ensure  that in terms of the empirical risk error, the TKRR and the  full KRR estimators have the same optimal convergence rate.  Finally, we provide the reader with some numerical simulations  that illustrate  the performance of our proposed TKRR estimator.\\
	

	
	\section{Introduction}
	
	For an integer $d\geq 1,$  let $\pmb X \subset \mathbb R^d$  be a compact set. Then, we consider the usual $d-$dimensional   non-parametric regression problem given by  
	\begin{equation}
		\label{model1}
		Y_i = f^*(X_i) + \varepsilon_i,\quad 1\leq i\leq n.
	\end{equation}
	Here,  the  $\pmb Y = [Y_i]^T_{1\leq i\leq n}$ is the vector of  real valued random responses,
	$f^*: \pmb X \rightarrow \mathbb R$ is the unknown  $d-$variate regression function to be estimated. The $ X_i \in \pmb X$ are the $\ n$ covariates.  The $ \varepsilon_i$ are   $n$ i.i.d. centered random variables with variance $\mathbb E\big[\varepsilon_i^2\big]=\sigma^2.$
	Also, we assume that the $\varepsilon_i$ are independent from the $X_i.$ In the literature, there is a variety of schemes providing different estimators for the regression function $f^*.$ The Kernel Ridge Regression (KRR) is among the popular schemes for providing fairly good estimate of $f^*$ under the hypothesis that this latter is well approximated by its projection over the associated reproducing kernel Hilbert space (RKHS). 
	
	Recently, there is a growing  interest in developing scalable methods  for  KRR algorithm with  the main objective: A drastic reduction of the computational complexity time of the KRR while maintaining the high convergence rate of this latter.    
	Among the existing scalable KRR estimators, we briefly describe the following four different improved KRR based estimator.    The first estimator is given by \cite{Amini}. It is briefly described as follows.
	The  covariates ${x_i}$ are assumed to be deterministic and   the regression function $f^*$ belongs to the RKHS $\mathcal H,$ associated with a positive semi-definite kernel $\mathbb K(\cdot,\cdot),$ defined on $\pmb X \times  \pmb X.$  Then, instead of using the full KRR estimator which involves the inversion of the regularized $n\times n$ Kernel Gram matrix, the author 
	has proposed an approximation scheme where the Gram matrix is  replaced by its best rank $r-$approximation, for some integer  $1\leq r\leq n.$ This is obtained by performing an SVD to the main Gram matrix and then considering only the main significant $r$ eigenvalues in the diagonal matrix and setting to  zero  the remaining $n-r$ eigenvalues.  The second scalable  KRR estimator is the Sketching based estimator given in \cite{Yang}, see also, \cite{Avron1}. More precisely, 
	given $n$ samples and a positive integer $m<n,$ the idea of the sketching is to replace the full kernel matrix by  $m-$dimensional randomized sketches of this latter. The projection dimension $m$ and the sketch matrix are  chosen in  such a way to maintain  the optimal convergence rate of the KRR estimator.   In particular,  the sketch dimension $m$ has to be proportional (modulo logarithmic factors) to the statistical dimension of the  kernel or Gram  matrix ${\bf K.}$ Recall that for a given regularization parameter $\lambda>0,$ the most 
	adopted definition of the statistical dimension of a kernel matrix ${\bf K}$ is given by 
	\begin{equation}\label{Eq1.1}
		d_{\bf K}(\lambda)= {\bf Tr }\big(({\bf K}+\lambda I_n)^{-1} {\bf K}\big),
	\end{equation}
	where ${\bf Tr}(A)$ denotes  the trace of a matrix $A.$ 
	The third scalable  KRR estimator is given in \cite{Rudi} and it is known as FALKON. This algorithm is based on the use of 
	random projections in order  to reduce the computational complexity and the  memory requirements of the conventional KRR. More precisely, given samples $\{x_i,y_i\}_{i=1}^n$, then only a subset $\{\tilde{x}_1\cdots\tilde{x}_M\}\subseteq\{x_1\cdots x_n\}$ with $M\ll n$ training points sampled uniformly are used to construct the estimator $$\widehat{f}_{\lambda,M}(x)=\sum_{i=1}^{M}\tilde{\alpha}_i \mathbb K(x_i,\tilde{x}_i).$$ Typically,  $M=O(\sqrt n)$ suffices for the FALKON algorithm to achieve  optimal statistical accuracy of the full $KRR$. The fourth but the pioneer of the scalable KRR is the Random Fourier Features (RFF), due to \cite{Rahimi}. The RFF requires that the positive-definite kernel $\mathbb K(\cdot,\cdot)$ is a convolution or stationary kernel of Bochner's type. Thanks to the structure of such a kernel, this latter is efficiently  approximated by its Monte-carlo approximation that involves $s$ Fourier features, 
	where $s\ll n.$ Consequently, the large size kernel matrix $K$ is substituted with an $s\times s$ positive semi-definite matrix involving $s$ Fourier features. 
	Recently, an extension of the RFF
	to  more general kernels has been provided in \cite{Zhu}.  A more detailed description of the RFF approach will be given in the next  section. \\
	
	We should mention  that our proposed TKRR and the RFF have some common points. Nonetheless, the TKRR does not require a convolution type structure of the positive-definite kernel. It is based on a substitution of the $n\times n$ random kernel   matrix 
	${\displaystyle B_n=\left[ \frac{1}{n} \mathbb K(X_i, X_j)\right]_{1\leq i,j\leq n}}$ by its  $n\times N$ main submatrix ${\displaystyle A_N=\left[ \frac{1}{n} \mathbb K(X_i, X_j)\right]_{\underset{1\leq j\leq N}{1\leq i\leq n}}},$ where $N\ll n.$ Moreover, the $n$ covariates $X_i$ are i.i.d. and drawn according to an unknown probability law $\rho(\cdot).$ We assume that  the measure $d\rho$ is absolutely continuous with respect the the usual probability measure $dP,$ associated with a Mercer's kernel $\mathbb K(\cdot,\cdot).$ More precisely, our TKRR is given as follows,
	\begin{equation}\label{Estimator0}
		\widehat{f}_{N,\lambda}(x)=\sum_{j=1}^{N}\widehat{\omega}_j \frac{1}{n} \mathbb K(X_j,x),\qquad 	\widehat{\pmb \omega} = \big[\widehat \omega_j\big]^T_{1\leq j\leq N} = \big(A_N^*A_N+\lambda I_N\big)^{-1} A_N^*  \pmb Y.
	\end{equation}

	In particular, we show that the choice of the optimal value of the truncation order is based on the behavior of the spectrum of the Hilbert-Schmidt integral operator $T_{\mathbb K},$ defined on $L^2(\pmb X, dP).$ The decay rate of the spectrum of $T_{\mathbb K},$ as well as the rank from which this decay holds true are crucial for the choice of optimal truncation order $N.$  
	We assume that for the  probability measure $dP$ on $\pmb X,$ the spectral properties of the 
	operator $T_{\mathbb K}$
	are known. The first main result of this work is the proof that 
	under some assumptions on the unknown sampling measure $d\rho(\cdot),$ the eigenvalues of full random kernel matrix, as well as the sequence of  singular values of the truncated $n\times N$ Gram matrix have 
	similar decay rate as the eigenvalues of the operator $T_{\mathbb K},$ when defined on $L^2(\pmb X, dP).$  The second main result of this work is an empirical $L^2-$risk error of  our TKRR estimator for solving  nonparametric (NP) regression problems.  In particular, we show that the proposed TKRR estimator has three main  advantages: Its empirical $L^2-$risk error is similar to the $L^2-$risk of the full KRR algorithm, it is much faster than the original KRR based algorithm and it is adapted for a whole class of unknown sampling probability laws. This makes our proposed distribution free  truncated KRR estimator  well adapted for analysis of driven data. \\
	
	This work is organized as follows. In section 2, we give some mathematical preliminaries on the classical KRR estimator as well as on the Random Fourier Features scalable KRR estimator. Then, we describe  our proposed TKRR estimator. In section 3, we prove a fairly  useful bound for   the eigenvalues of a compact  kernel integral operator, associated with an unknown probability measure, in terms of the eigenvalues of this operator when associated with its usual probability measure. Then, we use this result to estimate 
	the decay rate of the eigenvalues of a random kernel matrix, associated with a sampling set drawn from an unknown probability density.  In section 4, we give the empirical risk error of our TKRR estimator. Then, by using the spectral analysis of  section 3, we give  the convergence rate of the TKRR estimator under the assumptions that the spectrum of the associated kernel integral operator has an exponential or polynomial decay rates. In section 5, we give some spectral properties of the Sinc-kernel  and Gaussian kernel integral operators. Then, we provide the reader with   some numerical simulations  that illustrate  the theoretical properties  and  the performance of our proposed TKRR estimator.\\

	\section{Preliminaries on RKHS and related KRR based estimators}
	
	In this paragraph, we first give some preliminaries on the conventional KRR estimator. Then, we describe the scalable Random Fourier Features estimator (RFF). The proposed TKRR shares some ideas behind the RFF. Nonetheless, these two scalable KRR estimators have a fundamental difference. The RFF is based on a Monte-Carlo approximation of a convolution kernel of Bochner's type, while the TKRR is based on a truncation of a fairly general kernel random matrix.  Finally, we describe our proposed TKRR estimator. \\

	The KRR estimator for problem \eqref{model1} is briefly described as follows. We first recall that a  real-valued kernel $\mathbb K(\cdot,\cdot)$ defined on $\pmb X\times \pmb X$ is said to be a Mercer's kernel if it is continuous and positive semi-definite. For simplicity, we may assume in the sequel that $\mathbb K(\cdot,\cdot)$ is positive-definite. Moreover, let $ d\rho$ be a probability  measure with support  $\pmb X \subseteq \mathbb R^d$  and let $(\varphi_n)_{n\geq 1},  (\lambda_n)_{n\geq 1}$ be the orthonormal eigenfunctions and the corresponding  eigenvalues of the associated Hilbert-Schmidt  operator $T_{\mathbb K}.$ That is  $$T_{\mathbb K} \varphi_n(x) = \int_{\pmb X} \mathbb K(x,y)\varphi_n(y)\, d\rho(y)= \lambda_n \,\varphi_n(x),\quad x\in \pmb X.$$ Then, by Mercer's Theorem, we have: $ \mathbb K(x,y)=\sum_{n=0}^\infty \lambda_n \,\varphi_n(x) \varphi_n(y)$ for $x,y\in \pmb X$. 
	This  sum converges  uniformly over the compact set $\pmb X\times \pmb X.$ Moreover,  the associated RKHS  $\mathcal H$  is given by:
	\begin{equation}\label{RKHS}
		\mathcal H=\left\{ f\in L^2(\pmb X,d\rho),\,\, f=\sum_{n\geq 0} a_n(f) \varphi_n,\,\, \| f\|_{\mathcal H}^2=\sum_{n\geq 1} \frac{|a_n(f)|^2}{\lambda_n}<+\infty\right\}.
	\end{equation}
	For a given regularization parameter $\lambda >0,$ the KRR algorithm  consists in finding a solution $f_{\lambda}\in \mathcal H$ of the minimization problem
	\begin{equation}\label{Tikhonov1}
		\widehat f_K=\arg\min_{f\in \mathcal H} \left\{\frac{1}{n} \sum_{i=1}^n \Big(f(X_i)-Y_i)\Big)^2 +\lambda \|f\|^2_{\mathcal H}\right \}.
	\end{equation}
	  It is well known, see for example \cite{Smale1}, that a solution of \eqref{Tikhonov1} is given by the estimator
	\begin{equation}\label{Estimator3}
		\widehat f_K(x)= \sum_{i=1}^n \widehat c_{i} \mathbb K(X_i,x),
	\end{equation}
	where the expansion coefficients vector $\widehat{\mathbf C} = (\widehat c_i)_{1\leq i\leq n}$ is a solution of the system
	\begin{equation}\label{Tikhonov3}
		\left[ \Big[\mathbb K(x_i,x_j)\Big]_{1\leq i,j\leq n} + n \lambda I_n\right] \widehat{\mathbf C} =G_\lambda\, \widehat{\mathbf C} =\mathbf Y,\quad \mathbf Y= (Y_i)^T_{1\leq i\leq n}.
	\end{equation}
	Here, $I_n$ is the $n\times n$ identity matrix and $\frac{1}{n}[\mathbb K(x_i,x_j)]_{1\leq i,j\leq n}$ is  the random Gram matrix associated with the kernel $\mathbb K(\cdot,\cdot).$ This result is a consequence of the famous representer theorem. This  theorem is particularly useful in the sense that even if the RKHS associated with the minimization problem  \eqref{Tikhonov1} is of infinite dimension, the solution  \eqref{Estimator3} always lies in a finite dimensional space.
	Also, under the condition that the regularized random Gram matrix $G_\lambda$ given by \eqref{Tikhonov3} is invertible, the expansion coefficients vector $\widehat{\mathbf C}$ is given by 
	\begin{equation}\label{Tikhonov4}
		\widehat{\mathbf C}= G_\lambda^{-1} \mathbf Y = \left[ \Big[\mathbb K(x_i,x_j)\Big]_{1\leq i,j\leq n} + n\lambda I_n\right]^{-1} \mathbf Y,\quad  \mathbf Y= (Y_i)^T_{1\leq i\leq n}.
	\end{equation}
	
	As we have already mentioned, the previous conventional KRR estimator has  advantage of having a fast convergence rate. Nonetheless, it has the main drawback to require an intensive  computational load or a high time complexity. To overcome this drawback, some scalable KRR estimators have been recently proposed in the literature. The pioneer of these KRR approaches is the Random Fourier Features (RFF) estimator, proposed  in \cite{Rahimi}. This latter uses the convolution and  positive definite kernels of
	Bochner's type. That is those kernels of the form ${\displaystyle \mathbb K(x,y) = \kappa (x-y),\, x,y\in \pmb X,}$ where 
	$$\kappa(x-y) =\int_{\mathbb R^d} e^{-2i\pi  \omega^T\cdot (x-y)}\, dP(\omega)= \int_{\mathbb R^d} e^{-2i\pi  \omega^T \cdot x} e^{2i\pi  \omega^T \cdot y}\, dP(\omega).$$
	Here, $dP(\omega)$ is a probability measure. Typically, the kernel is real valued. Consequently, we may only consider the real part of this kernel. The idea behind the RFF is to replace the previous kernel by its following  Monte-Carlo approximation
	$$\widetilde{\mathbb K}(x,y)=\frac{1}{s} \sum_{j=1}^s Z^*_{\omega_j}(x)  Z_{\omega_j}(y)=Z_{\pmb \omega}^*(x) Z_{\pmb \omega}(y),\quad  Z_{\pmb \omega}(x)=\frac{1}{\sqrt{s}}
	\Big( \cos(2\pi  \omega_1^T \cdot x),\ldots, \cos(2\pi  \omega_s^T \cdot x)\Big)^T.$$
	Here, the $\omega_i$ are i.i.d. and drawn according to $P(\cdot)$ and $s$ is the number of Fourier features, with the usual assumption that $s\ll n.$  The $n\times n$ kernel matrix $\bf K$ and its Monte-Carlo  kernel matrix approximation $\widetilde{\bf K}$ are given by 
	\begin{equation}
		\label{KernelMatrix1}
		{\bf K}= \Big[\kappa(X_i-X_j)\Big]_{1\leq i,j\leq n},\quad  \widetilde{\bf K}= \pmb Z \pmb Z^* ,\quad \pmb Z=
		\Big[Z_{\pmb \omega}(X_1),\ldots,Z_{\pmb \omega}(X_n)\Big]^T\in \mathbb R^{n\times s}.
	\end{equation} 
	Note that since ${\displaystyle \mathbb E\big(\widetilde{\bf K}\big)= {\bf K},}$ then $\widetilde{\bf K}$ is an unbiased approximation of the kernel matrix ${\bf K}.$  The RFF estimator is given by 
	\begin{equation}
		\label{RFF_estimator}
		\widetilde f_{s,\lambda}(x)=\sum_{j=1}^s \widetilde \alpha_j Z_{\omega_j}(x)= Z^T_{\pmb \omega} \pmb{\widetilde \alpha},
		\qquad \pmb{\widetilde \alpha}=\big(\pmb Z^* \pmb Z+  s \lambda I_s\big)^{-1} \pmb{Z}^* \pmb Y.
	\end{equation} 
	That is $\pmb{\widetilde \alpha}\in \mathbb R^s$ is the solution of the minimization problem
	\begin{equation}\label{minimizationRFF}
		\pmb{\widetilde \alpha}=\argmin_{\alpha\in \mathbb R^s} 
		\frac{1}{n} \big\| \pmb Y- \pmb Z \alpha\big\|^2_{\ell_{2,n}}+s \lambda 
		\big\| \alpha\big\|^2_{\ell_{2,n}}.
	\end{equation} 
	Here, $\| \cdot\|_{\ell_{2,n}}$ is the usual Euclidean norm of $\mathbb R^n.$
	For more details, see for example \cite{Avron2}.  Recently, in \cite{Zhu}, an extension of the RFF has been provided and studied. It extends the RFF approach to those kernels of the form 
	\begin{equation}\label{generalRFF}
		\mathbb K(x,y)= \int_{\mathcal V} z(\omega, x) z(\omega, y)\, dP(\omega),
	\end{equation}
	where $z:{\mathcal V} \times \pmb X \rightarrow \mathbb  R $ is a continuous and bounded function.\\

	Next, we describe how to derive our proposed TKRR estimator. 
	For two positive integers $N\leq n$, we let $S_X=\{X_1,\ldots,X_n\}$ be a random sampling set following a probability law with a pdf $ P(\cdot),$ supported on $\pmb X.$ Let  $\mathbb K(\cdot,\cdot)$ be a continuous and positive-definite Mercer's kernel, so that the integral operator $T_{\mathbb K}$ defined on $L^2(\pmb X, dP)$ by 
	\begin{equation}\label{TK}
		T_{\mathbb K}(f)(x)=\int_{\pmb X} \mathbb K(x,y) f(y) d P(y)
	\end{equation}
	is a Hilbert-Schmidt operator with a countable set of positive eigenvalues, arranged in the decreasing order 
	$$\lambda_1 \geq \lambda_2\geq \cdots \geq \lambda_n \geq \cdots \geq 0. $$
	The eigenvalues $\lambda_n$ are associated to an orthonormal set of eigenfunctions $\{\varphi_n,\ n\in \mathbb N\}.$ Moreover, under the hypothesis that the integral operator $T_{\mathbb K}$ is one to one,  the $\varphi_n$ constitute an orthonormal basis of the Hilbert space $L^2(\pmb X, dP).$ 
	Recall that the Mercer's kernel has the spectral decomposition 
	\begin{equation}\label{Kernel1}
		\mathbb K(x,y)=\sum_{m=1}^\infty \lambda_m \varphi_m(x)\varphi_m(y),\qquad \forall\, x,y\in \pmb X.
	\end{equation}
	In the sequel, we  adopt the Notations:
	\begin{equation}\label{notation1}
		\pmb \varepsilon=(\varepsilon_1,\cdots,\varepsilon_n),\quad A_N^* A_N=K_N^2,\quad
		\| \pmb u  \|_n^2 = \frac{1}{n} \sum_{i=1}^n u_i^2,\quad \pmb u = (u_1,\ldots,u_n)^T\in \mathbb R^n.
	\end{equation}
	For   $\pmb Y\in \R^n,$ $ K_1\in \R^{n\times N},\, K_2\in \R^{N\times N}$ two positive definite matrices and $\lambda>0,$
	we aim to find 
	\begin{equation}\label{minmization1}
		\underset{\pmb \omega\in \R^N} \argmin \| K_1 \pmb \omega-\pmb Y\|_n^2+\lambda\omega^T \cdot  K_2\omega= \underset{\pmb \omega\in \R^N}  \argmin J_\lambda(\pmb \omega)
	\end{equation}
	Using the linearity of $K_1$ and the bi-linearity of $<\cdot,\cdot>,$ the usual inner product of $\mathbb R^n,$ together with the fact that $K_2^*=K_2,$ on gets 
	
	\begin{eqnarray*}
		\nabla J_\lambda(\pmb \omega)\cdot \pmb h &=&\frac{2}{n}<K_1\pmb \omega- \pmb Y,K_1 \pmb h>+2\lambda <K_2 \pmb \omega,\pmb h>
		=\frac{2}{n}<K_1^*K_1 \pmb \omega-K_1^* \pmb Y,\pmb h>+2\lambda<K_2 \pmb \omega,\pmb h>\\
		&=&2<\frac{1}{n}K_1^*K_1 \pmb \omega-\frac{1}{n}K_1^* \pmb Y+\lambda K_2\pmb \omega,\pmb h>
	\end{eqnarray*}
	So that ${\displaystyle \nabla J_\lambda(\pmb \omega)=0}$ whenever ${\displaystyle \Big(\frac{1}{n} K_1^*K_1+\lambda K_2\Big) \pmb \omega=\frac{1}{n}K_1^* \pmb Y.}$ That is 
	\begin{equation}\label{solutionmin1}
		\pmb \omega = \Big( K_1^*K_1+n \lambda K_2\Big)^{-1} K_1^* \pmb Y.
	\end{equation}
	The TKRR estimator is described as follows. For a given compact set $\pmb X\subset \mathbb R^d,$ a positive definite Mercer's kernel $\mathbb K$ defined on $\pmb X\times \pmb X$ and a set of $n$ i.i.d. sampling points following an unknown probability density supported on $\pmb X,$ we consider an integer $1\leq N\leq n$ and the rectangular random Gram matrix 
	\begin{equation}\label{GramMatrix1}
		K_1=A_N=\left[ \frac{1}{n} \mathbb K(X_i, X_j)\right]_{1\leq i\leq n, 1\leq j\leq N}
	\end{equation}
	Note  that $K_1^*K_1\in \R^{N\times N}$ is positive semi-definite and we can control the eigenvalues $\lambda_i\big(K_1^*K_1\big)$ of $K_1^*K_1$. Our TKRR estimator is denoted by $\widehat f_{N,\lambda}$ and it is given by the solution of the minimization problem \eqref{minmization1} with the previous Gram matrix $K_1=A_N$ and the matrix ${\displaystyle K_2=\frac{1}{n} I_N.}$ Hence, by using 
	\eqref{solutionmin1} with a convenient value for the  regularization parameter $\lambda>0,$ the estimator $\widehat f_{N,\lambda}$ is given by the following explicit formula
	\begin{equation}\label{Estimator1}
		\widehat{f}_{N,\lambda}(x)=\sum_{j=1}^{N}\widehat{\omega}_j \frac{1}{n} \mathbb K(X_j,x),\qquad 	\widehat{\pmb \omega} = \big(A_N^*A_N+\lambda I_N\big)^{-1} A_N^*  \pmb Y,\quad \pmb Y=\Big[Y_i\Big]^T_{1\leq i\leq n}.
	\end{equation}

	\section{Spectrum decay of  distributions free Gram matrices}
	
	In this paragraph, we first investigate the important problem of the  decay rate estimate of a random Gram matrix associated with a fairly general and unknown sampling probability law. Then, we extend this decay rate to the singular values of a rectangular random Gram matrix.  This last decay rate estimate will be used in the empirical risk error of our truncated KRR. More importantly, the issue of decay rate of the eigenvalues of Gram matrices is particularly useful in many kernels based algorithms for machine learning and data science applications. This is the case for example  for  kernels based unsupervised domain  learning problems, where  usually the source probability measure  is different from the target probability measure. This latter is  used to estimate the expected risk of the prediction of the output  $y\in \bf Y\subset \mathbb R$ from the input $\pmb x \in \mathbb R^d.$ For more details on this issue, the reader is referred to \cite{Gizewski} and the references therein. \\
	
	Typically, when a  KRR algorithm is used for solving a  learning problem, the associated Gram matrix is generally constructed from a set of $d-$dimensional sampling random vectors following an unknown probability measure $d\rho$ which is generally different from the usual   probability measure $dP,$ associated with the corresponding kenel integral operator.   For example, in the one dimensional case, for a positive parameter $c>0,$  the Sinc-kernel ${\displaystyle \mathbb K_c(x,y)=\frac{\sin(c(x-y))}{\pi (x-y)}},$ defined on $I=[-1,1]^2$  is associated to the uniform measure over $I,$ given by ${\displaystyle dP(x)= \frac{1}{2} \mathbf 1_{[-1,1]}(x) \, dx}.$ Also, for two positive real numbers $\xi, c >0,$ the Gaussian kernel ${\displaystyle \mathbb K_G(x,y)= e^{-\xi  |x-y|^2}}$ is associated to its usual  Gaussian  measure ${\displaystyle dP_c(x)= \sqrt{\frac{\pi}{2c}} e^{-2 c x^2}}$. For more details on the spectra and the super-exponential decay rate of the eigenvalues of  these two last kernels operators, the reader is referred to \cite{Bonami-Karoui1} and \cite{Ramussen}, respectively. \\
	
	It is well known that a Gram matrix has similar spectrum properties as the associated integral operator, see for example \cite{SCK1, Bonami-Karoui1}. Consequently, we are first interested in  comparing the spectra of two self-adjoint Hilbert-Schmidt  operators, $T_{{\mathbb K}}$ and $T_{\mathbb K,\rho},$  associated with the same continuous Mercer's kernel $\mathbb K(\cdot,\cdot)$ but acting on different Hilbert spaces $L^2(\pmb X, dP)$ and $L^2(\pmb X, d\rho).$
	Here, $dP$ is the original considered probability measure and $d\rho$ is an unknown probability measure supported on the compact set $\pmb X \subset \mathbb R^d.$ Moreover, we assume  the
	following two hypotheses:\\

	\noindent
	$ \bf H_0:$ The operator $T_{\mathbb K}: L^2(\pmb X, dP)\rightarrow  L^2(\pmb X, dP)$ is one-to-one.\\
	
	\noindent
	$\bf H_1:$    $d\rho$ is absolutely continuous with respect to $dP$ and 
		${\displaystyle \frac{d\rho}{dP}}$ is bounded on $\pmb X.$\\

	\noindent 
	We let $${\displaystyle \lambda_{1}\geq \lambda_{2}\geq \cdots \geq \lambda_{n}\geq \cdots\geq 0\mbox{\ and \ }  \lambda_{1,\rho}\geq \lambda_{2,\rho}\geq \cdots \geq \lambda_{n,\rho}\geq \cdots\geq 0}$$ denote  the  eigenvalues of $T_{\mathbb K}$ and $T_{{\mathbb K},\rho},$ arranged in the decreasing order. Also, we let $\big\{\varphi_n,\, n\in \mathbb N \big\}$ denote the set of the eigenfunctions of $T_{\mathbb K}.$ Note that by hypothesis $\bf H_0$ and the assumption that the kernel $\mathbb K(\cdot,\cdot)$ is a continuous Mercer's kernel, the $\varphi_n$ form an orthonormal basis of the Hilbert space $L^2(\pmb X, dP).$ In practice, regression estimators have to deal with $d-$dimensional  data samples. These samples are usually assumed to be i.i.d. and following  an unknown probability law with unknown pdf $\rho(\cdot).$ Given the previous continuous Mercer's kernel $\mathbb K(\cdot,\cdot),$ with associated self-adjoint Hilbert operator $T_{{\mathbb K}}$, having known spectral properties, we consider the following general case of $n\times n$ random Gram matrix
	\begin{equation}\label{GramM}
		B_n = \frac{1}{n} \Big[ \mathbb K(X_i,X_j) \Big]_{1\leq i,j\leq n},
	\end{equation}
	where the $X_i$ are i.i.d. random samples drawn according to  $\rho(\cdot).$ The following theorem provides us with interesting and useful  upper bounds for the unknown eigenvalues $\lambda_{k,\rho}$ as well as the expected value  for the tail of the trace of the random matrix $B_n.$
	
	\begin{theorem}\label{ThmCompare}
		Under the previous notations and hypotheses $\bf H_0$ and $\bf H_1,$ we have
		\begin{equation}\label{compare1}
			\lambda_{m,\rho} \leq \sum_{j\geq m} \lambda_{j} \|\varphi_j\|^2_\rho,\quad \forall\, m\geq 1. 
		\end{equation} 
		Moreover, for any integer $1\leq k\leq n,$ we have 
		\begin{equation}\label{compare2}
			\mathcal T_k(B_n)=\mathbb E_Z\Big[\sum_{j\geq k}  \lambda_j(B_n)\Big] \leq \sum_{j\geq k}\lambda_{j,\rho} \lesssim \lambda_{k} \|\varphi_k\|^2_\rho,\quad  \mbox{ whenever } \sum_{i\geq k} \frac{(i-k+1)\lambda_{i} \|\varphi_i\|^2_\rho}{\lambda_{k} \|\varphi_k\|^2_\rho} < +\infty.
		\end{equation} 
		Here, ${\displaystyle  \|\varphi_j\|^2_\rho=\int_{\pmb X} |\varphi^2_j(y)|\, d\rho(y).}$
	\end{theorem}
	
	\noindent
	{\bf Proof:} We first note that  hypothesis $\bf H_1$ and the well known and  powerful Radon-Nikodym theorem give us  the existence of the measure derivative ${\displaystyle \frac{ d\rho}{d P}.}$ Moreover, for any function $f\in  L^1(\pmb X, d\rho),$ we have the change of variable-type formula, 
	$$ \| f\|_\rho = \int_{\pmb X} |f(y)| d\rho(y) =\int_{\pmb X} |f(y)| \frac{d \rho(y)}{d P(y)}\, dP(y)=
	\int_{\pmb X} |f(y)| \frac{\rho(y)}{P(y)}\, dP(y).$$
	In particular, since the Mercer's kernel $\mathbb K(\cdot,\cdot)$ is continuous, then the eigenfunctions $\varphi_m$ of $T_{\mathbb K}$ are continuous over the compact set $\pmb X.$ Consequently, they are bounded and belong to $L^2(\pmb X, dP)\cap L^2(\pmb X, d\rho).$ That is the $\|\varphi_m\|_\rho$ are well defined for any integer $m\geq 1.$ Also, it is easy to see that  the function 
	${\displaystyle \psi_k: x\rightarrow \varphi_k(x) \sqrt{\frac{P(x)}{\rho(x)}}}$ belongs to $L^2(\pmb X, d\rho).$
	Next, for an integer $m\geq 1,$ consider the two  $m-$dimensional subspaces of $L^2(\pmb X, d\rho),$ given by 
	$$S_m=\mbox{Span} \Big\{\varphi_k(x) \sqrt{\frac{P(x)}{\rho(x)}};\, 1\leq k\leq m \Big\},\quad V_m=\mbox{Span} \Big\{\varphi_k(x); \, 1\leq k\leq m \Big\}.$$
	The linear independence of the $\varphi_k$ in $L^2(\pmb X, d\rho)$ is a consequence of their orthogonality in $L^2(\pmb X, dP)$ and the condition $\bf H_1.$ More precisely, if ${\displaystyle \sum_{k=1}^m a_k\varphi_k(x)=0,}$ then  ${\displaystyle \sum_{k=1}^m {a_k} \varphi_k(x) \frac{P(x)}{\rho(x)}=0.}$ Consequently for any $1\leq j \leq m,$ we have ${\displaystyle a_j= < \sum_{k=1}^m {a_k} \varphi_k(x) \frac{P(x)}{\rho(x)}, \varphi_j>_\rho=0.}$ That is the $\varphi_j$ are also linearly independent in $L^2(\pmb X, d\rho).$
	Next, we show that $S_m^{\perp}$ and $V_m^{\perp},$  the orthogonal of $S_m$ and $V_m$  in $L^2(\pmb X, d\rho)$ are  given by 
	$$S_m^{\perp}=\mbox{Span} \Big\{\varphi_k(x) \sqrt{\frac{P(x)}{\rho(x)}};\,  k\geq  m+1 \Big\},\quad V_m^{\perp}=\mbox{Span} \Big\{\varphi_k(x) \frac{P(x)}{\rho(x)};\,  k\geq  m+1 \Big\}.$$
	It is clear that from the change of variable-type formula, we have  
	$$\int_{\pmb X} \varphi_j(x) \varphi_k(x) \frac{P(x)}{\rho(x)} d\rho(x)= \int_{\pmb X} \varphi_j(x) \varphi_k(x) dP(x)=0,\quad \forall j\leq m,\, k\geq m+1.$$ That is any element from $S_m^{\perp}$ and $V_m^{\perp}$ is orthogonal to every element of $S_m$ and $V_m,$ respectively.  Next, assume that there exists $\psi \in L^2(\pmb X, d\rho)$ which is orthogonal to $S_m,$ then we check that $\psi\in S_m^{\perp}.$ Since the function $\psi(\cdot)\sqrt{\frac{\rho(\cdot)}{P(\cdot)}}$ belongs to $L^2(\pmb X, dP),$ then by writing for $1\leq k\leq m,$
	\begin{eqnarray*}
		0&=&\int_{\pmb X} \varphi_k(x)\sqrt{\frac{P(x)}{\rho(x)}} \psi(x)\, \rho(x) dx=\int_{\pmb X} \varphi_k(x) \psi(x) \sqrt{\frac{\rho(x)}{P(x)}}\, P(x) dx, 
	\end{eqnarray*}
	one concludes that $\psi(\cdot)\sqrt{\frac{\rho(\cdot)}{P(\cdot)}}$
	belongs to the orthogonal subspace $E_m^{\perp}=\mbox{Span} \Big\{\varphi_k(x), \, 1\leq k\leq m \Big\}^{\perp}.$ Moreover, since $\{\varphi_j,\, j\geq 1\}$ is an orthonormal basis of $L^2(\pmb X, dP),$ then $E_m^{\perp}=\mbox{Span}\{\varphi_k(x), \, k\geq  m+1 \Big\}.$ That is $\psi(\cdot)\in S_m^{\perp}.$ In a similar manner, we check that if $\psi\in L^2(\pmb X, d\rho)$ is orthogonal to $V_m,$ then $\psi\in \mbox{Span} \Big\{\varphi_k(x) \frac{P(x)}{\rho(x)};\, \,  k\geq m+1 \Big\}.$ This is a simple consequence of the fact that 
	\begin{equation}\label{orthogonality}
		0=\int_{\pmb X} \varphi_j(x) \psi(x) d\rho(x)=\int_{\pmb X} \varphi_j(x)\sqrt{\frac{P(x)}{\rho(x)}} \psi(x)\sqrt{\frac{\rho(x)}{P(x)}} d\rho(x),\quad 1\leq j\leq m
	\end{equation}
	and the second statement of hypothesis $\bf H_1$ that  implies that the function ${\displaystyle \psi(\cdot) \sqrt{\frac{\rho(\cdot)}{P(\cdot)}}\in L^2(\pmb X,d\rho).}$
	Consequently, from \eqref{orthogonality}, one concludes that 
	$\psi(\cdot)\sqrt{\frac{\rho(\cdot)}{P(\cdot)}}  \in S_n^{\perp}$ which implies that $\psi(\cdot) \in V_m^{\perp}.$
	
	Next, we use the Min-Max characterization of positive eigenvalues that are  arranged in decreasing order  of self-adjoint compact operator $T$ acting on a Hilbert space $\mathcal H.$ More precisely, 
	$$\lambda_m = \min_{U_{m-1}} \max_{f\in U_{m-1}^{\perp}, \|f\|_{\mathcal H}\leq 1} < Tf, f>_{\mathcal H}.$$
	Here, the  $U_m$ run over the $m-$dimensional subspaces of 
	$\mathcal H$ and $<\cdot,\cdot>_{\mathcal H}$ is the inner product of $\mathcal H.$ In particular, for the special case $\mathcal H=L^2(\pmb X, d\rho) $ and $U_{m-1}=V_{m-1},$ one gets
	$$\lambda_{m,\rho} \leq \max_{f\in V_{m-1}^{\perp}, \|f\|_\rho\leq 1} <T_{\mathbb{K},\rho}f,f>_\rho.$$ Since $f\in V_{m-1}^{\perp},$ then
	${\displaystyle f(x) =\sum_{k\geq m} \alpha_k \varphi_k(x) \frac{P(x)}{\rho(x)}.}$ Hence, by using the spectral decomposition of the Mercer's kernel $\mathbb K(\cdot,\cdot),$ together with the orthonormality of the $\varphi_k$ in $L^2(\pmb X, dP),$   one gets for $f\in V_{m-1}^{\perp},$
	\begin{eqnarray*}
		<T_{\mathbb K,\rho}f,f>&=&\int_{\pmb X} \left(\int_{\pmb X} \sum_{n\geq 1} \lambda_{n} \varphi_n(x)\varphi_n(y)\cdot \sum_{k\geq m} \alpha_k \varphi_k(y) \frac{P(y)}{\rho(y)} d\rho(y)\right) \sum_{j\geq m} \alpha_j \varphi_j(x) \frac{P(x)}{\rho(x)} d\rho(x)\\
		&=& \int_{\pmb X} \left(\sum_{n\geq 1,k\geq m} \varphi_n(x) \int_{\pmb X}  \lambda_{n} \cdot  \alpha_k \,  \varphi_n(y) \varphi_k(y) dP(y) \right) \sum_{j\geq m} \alpha_j \varphi_j(x) \frac{P(x)}{\rho(x)} d\rho(x)\\
		&=& \int_{\pmb X} \sum_{k\geq m} \lambda_{k} \cdot  \alpha_k \,  \varphi_k(x) \sum_{j\geq m} \alpha_j \varphi_j(x)  dP(x)\\
		&=& \sum_{k\geq m} \lambda_{k} \cdot \alpha_k^2.
	\end{eqnarray*}
	On the other hand, from the previous expansion of $f$ and Cauchy-Schwartz inequality,
	$$ |\alpha_k|^2 = |< f, \varphi_k>_\rho|^2 \leq \|f\|^2_\rho\,  \|\varphi_k\|^2_\rho.$$
	Consequently, 
	$$ \max_{f\in V_{m-1}^{\perp}, \|f\|_\rho\leq 1} <T_{\mathbb K,\rho}f,f>_\rho \leq \sum_{k\geq n} \lambda_{k} \cdot \|\varphi_k\|^2_\rho.$$
	This concludes the proof of \eqref{compare1}. To prove \eqref{compare2}, the first inequality can be found in \cite{SCK1},
	see also \cite{Bonami-Karoui1}. To prove the second inequality,  
	we use  \eqref{compare1} and get
	$$ \mathbb E_Z \Big[\sum_{j\geq k} \lambda_j(B_n)\Big] \leq \sum_{j\geq k} \sum_{i\geq j} \lambda_{i} \|\varphi_i\|^2_\rho \leq \sum_{i\geq k}  (i-k+1)\lambda_{i} \|\varphi_i\|^2_\rho \lesssim   \lambda_{k} \|\varphi_k\|^2_\rho.$$
	To check  the previous third inequality, it suffices to write
	$$\sum_{i\geq k} (i-k+1) \lambda_{i} \|\varphi_i\|^2_\rho= \lambda_{k} \|\varphi_k\|^2_\rho \cdot \sum_{i\geq k}  \frac{(i-k+1)\lambda_{i} \|\varphi_i\|^2_\rho}{\lambda_{k} \|\varphi_k\|^2_\rho} \lesssim 
	\lambda_{k} \|\varphi_k\|^2_\rho,$$
	whenever ${\displaystyle \sum_{i\geq k} \frac{(i-k+1)\lambda_{i} \|\varphi_i\|^2_\rho}{\lambda_{k} \|\varphi_k\|^2_\rho} < +\infty.}$ In a similar manner, one gets the previous second  inequality for the upper bound of $\mathbb E_Z \Big[\sum_{j\geq k} \lambda_j(B_n)\Big].$\\
	\qed
	
	We should mention that many  examples of Mercer's kernels used by KRR algorithms have sets of eigenvalues $\lambda_{m}$ that decay exponentially to zero. That is $\lambda_{m} \lesssim e^{-bm},$ for some $b>0.$ The magnitudes $\|\varphi_m\|_\infty$ of their associated set of eigenfunctions $\varphi_m$ are bounded by a quantity of $O(m^{a/2}),$ for some $a\geq 0.$ For the unknown probability measure $d\rho$ which we assume to be  absolutely continuous with respect to the original measure $dP,$ and according to \eqref{compare1}, we have 
	for any $m\geq 0,$
	\begin{equation}\label{bound1}
		\lambda_{m+1,\rho}\lesssim \sum_{k\geq m+1} k^a e^{-b k} \lesssim \int_m^{\infty} t^a e^{- b t}\, dx =\frac{1}{b^{a+1}} \int_{b m}^{\infty} u^{a} e^{-u} \, du =\frac{1}{b^{a+1}} \Gamma(a+1,bm).
	\end{equation}
	Here, $\Gamma(a+1,bm)$ is the upper incomplete Gamma function. It is well known that for $x > a-1,$ we have ${\displaystyle \Gamma(a,x) \leq \frac{x^a e^{-x}}{x+1-a}.}$ Consequently, the eigenvalues of such an  integral operator with a given Mercer's kernel and an eventually unknwon probability measure $d\rho,$ have the following upper bounds,
	$$\lambda_{m+1,\rho} \lesssim  \frac{ (bm)^{a+1} e^{-bm}}{bm-a}\lesssim m^a e^{-b m}.$$
	Moreover, in this case for a  Gram matrix $B_n,$ given by \eqref{GramM} and according to  \eqref{compare2}, the expected value of its trace tail is bounded by 
	\begin{equation}\label{bound2} 
		\mathbb E_Z\Big[\sum_{j\geq m}  \lambda_j(B_n)\Big] \lesssim 
		m^{a} e^{- b m}.
	\end{equation}
	
	Next, for a positive integer $N\leq m,$ we extend the decay rate estimate \eqref{compare2} for the eigenvalues of the Gram matrix $B_n$ to the case of the singular values of $A_N,$  the principal  rectangular $m\times N$ sub-matrix  of $B_n.$ That is 
	\begin{equation}\label{GramM2}
		A_N = \frac{1}{n} \Big[ \mathbb K(X_i,X_j) \Big]_{1\leq i\leq n,\, 1\leq j\leq N}.
	\end{equation}
	We recall that  the $X_i$ are i.i.d. random samples following an unknown probability law with  pdf $\rho(\cdot).$ For this purpose, we use the following interlacing inequalities for the singular values of sub-matrices, given in \cite{Thompson}.\\
	
	\noindent
	{\bf Interlacing inequalities of singular values}\, \cite{Thompson}: {\it Let $B$ be an $m\times n$ matrix with singular values 
		$$\sigma_1(B)\geq \sigma_2(B)\geq \cdots\geq \sigma_{\min(m,n)}(B).$$ Let $D$ be a $p\times q$ submatrix of $B,$ with singular values 
		$$\mu_1(D)\geq \mu_2(D)\geq \cdots\geq \mu_{\min(p,q)}(D).$$ Then, we have 
		\begin{equation}\label{Interlacing}
			\sigma_i(B)\geq \mu_i(D)\geq \sigma_{i+(m-p)+(n-q)}(B),\qquad \forall\, 1\leq i\leq \min(p+q-m,p+q-n).
		\end{equation}
	}
	By using the first inequality of \eqref{Interlacing} and the inequalities \eqref{compare2} of Theorem 2, one gets the following Corollary.
	
	\begin{corollary} Let $B_n$ and $A_N$ be the Gram matrix and its principal sub-matrix, given by \eqref{GramM} and \eqref{GramM2}, respectively. Under the hypotheses of Theorem 2 and the previous notation,  if ${\displaystyle \kappa_1=\sup_{x\in \pmb X} \mathbb K(x,x),}$ then for any $1\leq m\leq N,$ we have 
		\begin{equation}\label{SingularValues1}
			\mathcal T_M=\mathbb E_Z\Big[\sum_{j\geq m} \mu_j(A_N) \Big]   \lesssim  \lambda_m \|\varphi_m\|^2_\rho
		\end{equation}
		and
		\begin{equation}\label{SingularValues2}
			\mathbb E_Z\Big[\sum_{j\geq m} \mu_j^2(A_N) \Big]  \lesssim \frac{\kappa_1}{m} \lambda_m \|\varphi_m\|^2_\rho.
		\end{equation}
	\end{corollary}
	
	\noindent
	{\bf Proof:} By using the interlacing property \eqref{Interlacing} and since the Gram matrix $B_n$ is positive  semi-definite, that is  $\sigma_i(B_n)=\lambda_i(B_n)$, then we have 
	\begin{equation}\label{SingularValues}
		\mu_i(A_N) \leq \lambda_i(B_n),\qquad \forall\, 1\leq i\leq N.
	\end{equation}
	By using the previous inequality and \eqref{compare2}, one gets 
	$$\mathbb E_Z\Big[ \sum_{j\geq m} \lambda_j(B_n) \Big]\leq    \sum_{j\geq m} \lambda_{j,\rho} \lesssim  \lambda_m \|\varphi_m\|^2_\rho.$$
	To check \eqref{SingularValues2}, we note that
	since ${\displaystyle \lambda_1(B_n)\geq \lambda_2(B_n)\geq \cdots \geq \lambda_n(B_n)\geq 0 }$ and since 
	$$ {\bf Tr}(B_n)= \sum_{i=1}^n \lambda_i(B_n)= \frac{1}{n} \sum_{i=1}^n \mathbb K(Z_i,Z_i) \leq \kappa_1,$$
	then $m \cdot  \lambda_m(B_n) \leq \kappa_1.$ Hence, by using \eqref{SingularValues}, one gets ${\displaystyle \mu_m(A_N)\leq
		\frac{\kappa_1}{m}}$. Consequently, one gets 
	\begin{eqnarray*}
		\mathbb E_Z\Big[\sum_{j\geq m} \mu_j^2(A_N) \Big] \leq
		\mathbb E_Z\Big[ \lambda_m(B_n)\sum_{j\geq m} \mu_j(A_N) \Big]
		\leq  \frac{\kappa_1}{m} \mathbb E_Z\Big[ \sum_{j\geq m} \lambda_j(B_n) \Big] \lesssim \frac{\kappa_1}{m} \lambda_m \|\varphi_m\|^2_\rho.
	\end{eqnarray*}

	\section{Empirical risk error analysis and convergence rate of the TKRR estimator}
	
	In this paragraph, we study the empirical risk error of our TKRR estimator $\widehat f_{N,\lambda},$ given by \eqref{Estimator1}
	and  under the usual assumption that 
	the true regression function $f^*$ belongs to $\mathcal H,$ the RKHS associated  with the positive definite Mercer's kernel $\mathbb K(\cdot,\cdot).$ Recall  that this empirical  risk error is given by 
	\begin{equation}
		\mathcal R(\widehat f_{N,\lambda})= \mathbb E\Big( \big\|\big(\widehat f_{N,\lambda}(X_i)-f^*(X_i)\big)_i\big\|_n^2\Big)=\mathbb E\Big(\frac{1}{n} \sum_{i=1}^n \big(\widehat f_{N,\lambda}(X_i)-f^*(X_i)\big)^2\Big).
	\end{equation}
	Here, $\mathbb E = \mathbb E_{X\times \varepsilon} =\mathbb E_{X} \mathbb E_{\varepsilon},$ the product of the expectations with respect to the independent random variables $X_i$ and $\varepsilon_i,$ respectively. Recall that the given i.i.d. random sampling  points  $X_i$ are drawn according to  $\rho(\cdot).$ We consider the following finite dimensional subspace of the RKHS $\mathcal H,$ given by 
	\begin{equation}\label{subspaces}
		\mathcal H_N =\mbox{Span}\Big\{ \frac{1}{n} \mathbb K(X_j,\cdot),\, 1\leq j\leq N \Big\},\qquad \mathcal F_N=\mbox{Span}\Big\{ \varphi_j(\cdot),\,   1\leq j\leq N \Big\}.
	\end{equation}
	We assume that $\mathcal H$ is a subspace of $L^2(\pmb X, dP),$ 
	and the true regression function $f^*\in \mathcal H.$ That is 
	${\displaystyle f^*=\sum_{j=1}^\infty a_j(f^*) \varphi_j,\,\, a_j(f^*)=<f^*,\varphi_j>_{\mathcal H}.}$ Also, we consider the following orthogonal projection of $f^*$ over the subspaces $\mathcal H_N$ and $\mathcal F_N,$ given by
	\begin{equation}\label{projections}
		\widetilde \pi_N f^*(\cdot) = \sum_{j=1}^N \widetilde \omega_j \frac{1}{n} \mathbb K(X_j,\cdot),\qquad \pi_N f^*(\cdot) = \sum_{j=1}^N a_j(f^*) \varphi_j(\cdot).
	\end{equation}
	Recall that $\widetilde \pi_N f^*$ and $\pi_N f^*$ are the best approximations of $f^*$ in the $L^2(dP)-$norm, by elements from $\mathcal H_N$ and $\mathcal F_N,$ respectively.
	In the sequel, we use the following two technical hypothesis:\\
	
	\noindent
	${\bf H_2:}$ The random projection matrix ${\displaystyle \Big[\varphi_j(X_i)\Big]_{1\leq i\leq n, 1\leq j\leq N}}$ has full rank.\\
	Next, we consider the truncated projection kernel,
	\begin{equation}\label{KernelN}
		\mathbb K_N(x,y)=\sum_{m=1}^N \lambda_m \varphi_m(x) \varphi_m(y),\quad x,y\in \pmb X.
	\end{equation}
	It is easy to see that under condition $\bf H_2,$ the projection $\pi_N f^*$ has the following second equivalent expansion form
	\begin{equation}\label{projectionKN}
		\pi_N f^*(x)= \sum_{j=1}^N \omega'_j \frac{1}{n}\mathbb K_N(X_j,x),\quad x\in \pmb X.
	\end{equation}
	Finally, to study the empirical risk error of our TKRR estimator, we  also need the following hypothesis on the coefficients of the projections $\widetilde \pi_N f^*$ and  $\pi_N f^*,$ given by \eqref{projections} and 
	\eqref{projectionKN}:\\
	
	\noindent
	${\bf H_3:}$ ${\displaystyle \max\big(\big\| \widetilde{\pmb \omega}\big\|_n^2, \big\| {\pmb \omega'}\big\|_n^2\big)\leq 1,\quad}$ where ${\displaystyle
		\quad \big\| \widetilde{\pmb \omega}\big\|_n^2=\frac{1}{n} \sum_{j=1}^N \widetilde \omega_j^2,\quad \big\| {\pmb \omega'}\big\|_n^2 = \frac{1}{n} \sum_{j=1}^N  {\omega'}_j^2.}$\\
	Note that the previous hypothesis can be relaxed by substituting the upper bound $1$ by any  constant $C>0.$
	The following technical lemma will be needed in the proof of our TKRR empirical risk error. It shows that $\widetilde \pi_N f^*$ is well approximated by $\pi_N f^*,$ whenever $f^*\in \mathcal H.$
	
	\begin{lemma}
		Under the previous notation and hypotheses ${\bf H_1}, \, {\bf H_2}$ and ${\bf H_3},$ let $f\in \mathcal H,$ then we have 
		\begin{equation}\label{projections2}
			\mathbb E\Big(\big\| \big((\pi_N f - \widetilde \pi_N f)(X_i)\big)_i\big\|^2_n\Big) \lesssim \frac{1}{N}\sum_{k\geq N+1} \lambda_k \|\varphi_k\|^2_\infty + \lambda_{N+1} \|f\|^2_{\mathcal H}. 
		\end{equation}
	\end{lemma}
	
	\noindent
	{\bf Proof:} It is well known that if $f\in \mathcal H,$ then we have 
	\begin{eqnarray*}
		\|f-\pi_N f\|^2_{P}&=& \sum_{k\geq N+1} |a_k(f)|^2 = \sum_{k\geq N+1} \lambda_k \frac{|a_k(f)|^2}{\lambda_k}\\
		&\leq& \lambda_{N+1} \sum_{k\geq N+1}  \frac{|a_k(f)|^2}{\lambda_k} \leq \lambda_{N+1} \|f\|_{\mathcal H}^2.
	\end{eqnarray*}
	Here, ${\displaystyle a_k(f)= <f,\varphi_k>_{P},}$ where $<\cdot,\cdot>_P$ is the usual inner product of $L^2(\pmb X, dP).$ Also, by using $\bf H_1,$ we have for any $g\in \mathcal H,$ 
	\begin{equation}
		\| g\|^2_\rho= \int_{\pmb X} |g(x)|^2\, d\rho(x) \lesssim \int_{\pmb X} |g(x)|^2\, dP(x)=\| g\|^2_{L^2(dP)}.
	\end{equation}
	On the other hand, for $1\leq i\leq n,$ we have 
	$$\mathbb E\Big(\big((\pi_N f - \widetilde \pi_N f)(X_i)\big)^2\Big)=\int_{\pmb X} |\widetilde \pi_N(f)(x)-\pi_N f(x)|^2 \, d\rho(x) \lesssim \| \widetilde \pi_N(f) -\pi_N(f)\|_{L^2(dP)}^2.$$
	Moreover, since $\widetilde \pi_N$  is a projection operator, then 
	we have 
	\begin{equation}
		\big\| \widetilde \pi_N \big(f-\pi_N f\big)\big\|_{L^2(dP)} \leq \|f-\pi_N f\|_{L^2(dP)} \leq \sqrt{\lambda_{N+1}} \|f\|_{\mathcal H}.
	\end{equation}
	Also, since $\widetilde \pi_N( \pi_N f)$ is the best approximation of $\pi_N f$ by elements from the subspace $\mathcal H_N$ and since 
	the non truncated kernel  expansion 
	${\displaystyle \sum_{j=1}^N \omega'_j \frac{1}{n}\mathbb K(X_j,x) \in \mathcal H_N,}$ then we have 
	\begin{eqnarray*}
		\big\|\widetilde \pi_N( \pi_N f)-\pi_N f\big\|_{L^2(dP)}^2 &\leq & \frac{1}{n^2} \Big\| \sum_{j=1}^N \omega'_j \sum_{k\geq N+1} \lambda_k \varphi_k (X_j) \varphi_k(\cdot)\Big\|_{L^2(dP)}^2
		\leq  \Big(\frac{1}{n}\sum_{k\geq N+1} \lambda_k \|\varphi_k\|^2_\infty\Big) \|\pmb \omega'\|_n^2,
	\end{eqnarray*}
	Finally, since 
	$$\| \widetilde \pi_N f -\pi_N f\|_{L^2(dP)}\leq \|\widetilde \pi_N (f-\pi_N f)\|_{L^2(dP)}+ \|\widetilde \pi_N(\pi_N f)-\pi_N f\|_{L^2(dP)},$$
	then, by using the previous inequalities, one gets the desired result 
	\eqref{projections2}.\\
	
	The following theorem provides us with an upper bound  for the empirical risk error of our TKRR estimator. 
	
	\begin{theorem}
		Let $n,N>0$ be two positive integers and $\lambda>0.$ Under the previous notation,  and hypotheses ${\bf H_1}, \, {\bf H_2}$ and ${\bf H_3},$  we have for $f^*\in \mathcal H,$
		\begin{eqnarray}\label{Riskerror}
			\mathbb{E}\Big( \big\|\big(\widehat{f}_{N,\lambda}-f^*\big)(X_i)\big\|_n^2\Big)&\leq&\frac{\lambda}{2}+2\frac{\sigma^2}{n}\sum_{j=1}^N\mathbb E_X\Big[\Big(\frac{\mu_j^2}{\mu_j^2+\lambda}\Big)^2\Big]\nonumber\\
			&&\qquad\qquad\qquad +\, C_1\, \Big(\lambda_{N+1} \|f^*\|^2_{\mathcal H}+\frac{1}{N}\sum_{k\geq N+1} \lambda_k \|\varphi_k\|^2_\infty  \Big).
		\end{eqnarray}
		Here, $C_1>0$ is a uniform positive constant.
	\end{theorem}
	\noindent
	\textbf{Proof :} Let $\widetilde \pi_N$ and $\pi_N$ be the projection operators defined by \eqref{projections}, then we have
	\begin{equation}
		\label{Eqq3.1}
		f^*(X_j)=\widetilde \pi_N f^*(X_j) +\Big( f^*(X_j)-\pi_N f^*(X_j)\Big)+\Big(\pi_N f^*(X_j)-\widetilde \pi_N f^*(X_j)\Big),\quad 1\leq j\leq n.
	\end{equation}
	Hence, from \eqref{Estimator1}, we have
	\begin{equation}\label{Eqq3.2}
		\Big[\widehat f_{N,\lambda}(X_i)\Big]^T_{1\leq i\leq n} =  A_N \big( A_N^* A_N+\lambda I_N\big)^{-1} A_N^* \Big[ f^*(X_j)+\varepsilon_j\Big]_{1\leq j\leq n}^T.
	\end{equation}
	By combining the previous two identities, one gets with ${\displaystyle F_N= \big(\widetilde \pi_N f^*(X_i)\big)_{1\leq i\leq n},}$
	\begin{eqnarray}\label{Eqq3.3}
		\lefteqn{\Big[\widehat f_{N,\lambda}(X_i)\Big]^T_{1\leq i\leq n} =  A_N \big(A_N^* A_N+\lambda I_N\big)^{-1} A_N^* \Big[ \big(\widetilde \pi_N f^*(X_j) +\varepsilon_j)+}\nonumber\\
		&& \qquad\qquad\qquad\qquad\qquad\qquad  +\big( f^*(X_j)-\pi_N f^*(X_j)\big)+\big(\pi_N f^*(X_j)-\widetilde \pi_N f^*(X_j)\big)\Big]_{1\leq j\leq n}^T.
	\end{eqnarray}
	Let $\|A\|_2$ denote the $2-$norm of a  matrix $A.$ By using \eqref{Eqq3.1} and \eqref{Eqq3.3}, together with the fact that $\big\|A_N \big(A_N^* A_N+\lambda I_N\big)^{-1}A_N^*\big\|_2\leq 1,$ one gets
	\begin{eqnarray}\label{Eqq3.4}
		\big\|\big(\widehat f_{N,\lambda}(X_i)-f^*(X_i)\big)_i \big\|_n&\leq& \big\|A_N(K_N^2+\lambda I_N)^{-1}A_N^* (F_N+\pmb \varepsilon)^T- F_N^T\|_n\nonumber\\
		&+&2\Big(\big\|\big(f^*(X_i)-\pi_N f(X_i)\big)_i\big\|_n+ \big\|\big(\pi_N f^*(X_i)-\widetilde\pi_N f(X_i)\big)_i\big\|_n\Big).
	\end{eqnarray}
Since $f^*\in \mathcal H,$ then as in the proof of  lemma~1, we have 
	\begin{eqnarray}\label{Eqq3.5}
		\mathbb E\left(\big\|\big(f^*(X_i)-\pi_N f(X_i)\big)_i\big\|^2_n\right) &=& \frac{1}{n} \sum_{i=1}^n \int_{\pmb X} |f^*(x)-\pi_N f^*(x)|^2 d \rho(x)\nonumber\\
		& \lesssim & \int_{\pmb X} |f^*(x)-\pi_N f^*(x)|^2 d P(x)= \| f^*-\pi_N f^*\|_{L^2(dP)}^2\nonumber\\
		& \lesssim &\lambda_{N+1} \| f^*\|^2_{\mathcal H}.
	\end{eqnarray}
	Moreover, from Lemma~1, we have 
	\begin{equation}\label{Eqq3.6}
		\mathbb E\Big(\big\| \big((\pi_N f^* - \widetilde \pi_N f^*)(X_i)\big)_i\big\|^2_n\Big) \lesssim \frac{1}{N}\sum_{k\geq N+1} \lambda_k \|\varphi_k\|^2_\infty + \lambda_{N+1} \|f^*\|^2_{\mathcal H}. 
	\end{equation}
	Next, we estimate the following main quantity in the empirical risk error, given by
	$$R_{n,N}=\mathbb E\left(\big\| \big(A_N \big( A_{N}^* A_{N}+\lambda I_N\big)^{-1} A_{N}^* \big[ \big(\widetilde \pi_N f^*(X_j)+\varepsilon_j\big)_j \big]-\widetilde \pi_N f^*(X_i)\big)_i\big\|^2_n\right).$$
	By using some techniques similar to  those used in the proof of the main theorem of \cite{Amini}, but in the present random sampling setting,  one has   the following singular values decomposition of $A_N=U\Sigma_{n,N}V^*$ where  $U, V$ are orthogonal matrices and   $\Sigma_{n,N}$ is an $n\times N$ rectangular diagonal matrix. Consequently, we have 
	$$K_N^2=A_N^* A_N=\big(U\Sigma_{n,N}V^*\big)^*\big(U\Sigma_{n,N}V^*\big)=V\Sigma_{n,N}^*\Sigma_{n,N}V^*.$$
	Let $D_N=\Sigma_{n,N}^*\Sigma_{n,N}.$ Then,
	$$(K_N^2+\lambda I_N)^{-1}=\big(V(D_N+\lambda I_N)\big)V^*)^{-1}=V(D_N+\lambda I_N)^{-1}V^*$$
	and 
	\begin{eqnarray*}
		A_N(K_N^2+\lambda I_N)^{-1}A_N^*&=&U\Sigma_{n,N}V^*V(D_N+\lambda I_N)^{-1}V^*V\Sigma_{n,N}^*U^*\\
		&=&U\Sigma_{n,N}(D_N+\lambda I_N)^{-1}\Sigma_{n,N}^*U^*.
	\end{eqnarray*}
	So that   
	\begin{eqnarray}\label{SVD}
		\big\|A_N(K_N^2+\lambda I_N)^{-1}A_N^* (F_N+\pmb \varepsilon)^T- F_N^T\|_n&\leq& \|U\Sigma_{n,N}(D_N+\lambda I_N)^{-1}\Sigma_{n,N}^*U^* F_N^T-F_N^T\|_n \nonumber\\
		&&+\|U\Sigma_{n,N}(D_N+\lambda I_N)^{-1}\Sigma_{n,N}^* U^* \pmb \varepsilon^T\|_n.
	\end{eqnarray}
	Next let $\|\cdot\|_{\ell_2,n}$ and $\|\cdot\|_2$ denote the usual $2-$norm of $\mathbb R^n$  and a rectangular matrix, respectively. In particular, for an $n\times N$ matrix $A$, $\|A\|_2=\sigma_1(A),$ the largest singular value of $A.$ Also, since
	${\displaystyle \widetilde \pi_N f^*(X_i)= \sum_{j=1}^N \widetilde \omega_j\, \frac{1}{n} \mathbb K(X_i,X_j),\qquad 1\leq i\leq n,}$
	then we have 
	$$ F_N^T =\Big(\widetilde \pi_N f^*(X_i)\Big)^T_{1\leq i\leq n}= A_N\,  \widetilde{\pmb{\omega}} = U\Sigma_{n,N}V^* \,\widetilde{\pmb{\omega}} ,\qquad \widetilde{\pmb{\omega}}=(\widetilde\omega_1,\ldots,\widetilde\omega_N)^T.$$
	Moreover, since the $2-$norm is invariant under  multiplication
	by an orthogonal matrix, one gets  
	\begin{eqnarray}\label{Eqq3.8}
		\big\|U\Sigma_{n,N}(D_N+\lambda I_N)^{-1}\Sigma_{n,N}^*U^* F_N^T-F_N^T\big\|_n^2
		&=&
		\frac{1}{n}\big\|U\big(\Sigma_{n,N}(D_N+\lambda I_N)^{-1}\Sigma_{n,N}^*-I_n\big)\Sigma_{n,N}V^*\widetilde{\pmb{\omega}}\big\|_{\ell_2,n}^2 \nonumber\\ 
		&=&\frac{1}{n}\big\|\big(\Sigma_{n,N}(D_N+\lambda I_N)^{-1}\Sigma_{n,N}^*-I_n\big)\Sigma_{n,N}V^*\widetilde{\pmb{\omega}}\big\|_{\ell_2,n}^2\nonumber\\
		&\leq &\big\|\big(\Sigma_{n,N}(D_N+\lambda I_N)^{-1}\Sigma_{n,N}^*-I_n\big)\Sigma_{n,N}\big\|_2^2 \cdot \frac{1}{n} \|V^*\widetilde{\pmb{\omega}}\|_{\ell_2,n}^2\nonumber\\
		&\leq &\big\|\big(\Sigma_{n,N}(D_N+\lambda I_N)^{-1}\Sigma_{n,N}^*-I_n\big)\Sigma_{n,N}\big\|_2^2 \cdot  \|\widetilde{\pmb{\omega}}\|_n^2.
	\end{eqnarray}
	Let ${\displaystyle \phi_{n}^\lambda=\Sigma_{n,N}(D_N+\lambda I_N)^{-1}\Sigma_{n,N}^*-I_n,}$
	then we  have
	\begin{equation}
		\label{Eqq3.9}
		\big\|\big(\Sigma_{n,N}(D_N+\lambda I_N)^{-1}\Sigma_{n,N}^*-I_n\big)\Sigma_{n,N}\big\|_2^2=\big|\sigma_1\big(\phi_{n}^\lambda\Sigma_{n,N}\big)\big|^2.
	\end{equation}
	On the other hand, we have  ${\displaystyle \Sigma_{n,N}=\big[s_{i,j}\big]_{\underset{1\leq j\leq N}{1\leq i\leq n}}, \quad s_{i,j}=\mu_j\delta_{i,j},}$ where the $\mu_j=\sqrt{\lambda_{j}(A_N^* A_N)}$ are the singular values of $A_N.$  Consequently, we have 
	$$\phi_{n,N}^\lambda\Sigma_{n,N}=\big[\gamma_{i,j}\big]_{\underset{1\leq j\leq N}{1\leq i\leq n}},\quad \gamma_{i,j}=\Big(\frac{{\mu}_j^2}{{\mu}_j^2+\lambda}-1\Big)\mu_j\delta_{i,j}=\frac{-\lambda \mu_j}{{\mu}_j^2+\lambda}\delta_{i,j}.$$
	Hence, 
	${\displaystyle |\sigma_1\big(\phi_{n}^\lambda\Sigma_{n,N}\big)|^2=\underset{1\leq j\leq N}{\max}\big(\frac{\lambda}{{\mu}_j^2+\lambda}\big)^2\mu_j^2.}$
	Since  ${\displaystyle \underset{x\geq 0}{\sup }\frac{\lambda x}{x^2+\lambda}=\frac{\sqrt{\lambda}}{2},}$ then ${\displaystyle \underset{1\leq j\leq N}{\max}\Big(\frac{\lambda}{{\mu}_j^2+\lambda}\Big)^2\mu_j^2=\frac{\lambda}{4}.}$ That is 
	\begin{equation}\label{Eqq3.10}
		\big|\sigma_1\big(\phi_{n}^\lambda\Sigma_{n,N}\big)\big|^2\leq \frac{\lambda}{4}.
	\end{equation}
	Hence, by using \eqref{Eqq3.8}--\eqref{Eqq3.10}, together with hypothesis
	${\bf H_3},$ one gets 
	\begin{equation}\label{Eqq3.11}
		\mathbb E\left(\big\|U\Sigma_{n,N}(D_N+\lambda I_N)^{-1}\Sigma_{n,N}^*U^* F_N^T-F_N^T\big\|_n^2\right) \leq \frac{\lambda}{4}.
	\end{equation}
	Next, to bound the expectation of the second term of the right hand-side
	of \eqref{SVD}, we first note that if $B_n=\big[b_{ij}\big]_{1\leq i,j\leq n}$ is an $n\times n$ matrix with real coefficients, then its Frobenius  norm $\|B_n\|_F$ is given by 
	$$ \|B_n\|^2_F=\sum_{i,j=1}^n (b_{ij})^2= \mbox{{\bf Tr}} (B_n^* B_n)=\sum_{i=1}^n \sigma_i^2(B_n),$$
	where, the $\sigma_i(B_n)$ are the singular values of $B_n.$ Moreover, the Frobenius norm is invariant under left or right multiplication by orthogonal matrices. We consider the special choice of 
	\begin{equation}
		\label{EQQ3.11}
		B_n= U\Sigma_{n,N}(D_N+\lambda I_N)^{-1}\Sigma_{n,N}^* U^*.
	\end{equation}
	In this case, $B_n$ has at most rank $N$ and 
	$$ \|B_n\|_F^2= \sum_{i,j=1}^n b^2_{ij}= \sum_{j=1}^N \Big(\frac{\mu_j^2}{\mu_j^2+\lambda}\Big)^2.$$
	Also, by using the fact that the $\varepsilon_i$ are   i.i.d. and  independent from the $X_i$ with 
	$\mathbb E_{\varepsilon}(\varepsilon_i)=0$ and  $\mathbb E_{\varepsilon}(\varepsilon^2_i)=\sigma^2,$
	one can easily check that 
	\begin{eqnarray}\label{Eqq3.12}
		\frac{1}{n}\mathbb E_{X\times\varepsilon} \left(\big\|B_n \pmb{\varepsilon}^T\big\|^2_{2}\right)&=&\frac{1}{n} \mathbb E_X \mathbb E_{\varepsilon} \sum_{i=1}^n \left(\sum_{j=1}^n b_{ij} \varepsilon_j\right)^2 \nonumber\\
		&=& \frac{1}{n} \mathbb E_X \left[\sum_{i,j=1}^n  b^2_{ij} \mathbb E_{\varepsilon}(\varepsilon_j^2)\right]=\frac{\sigma^2}{n}\mathbb E_X \left[\sum_{j=1}^N \Big(\frac{\mu_j^2}{\mu_j^2+\lambda}\Big)^2\right]. 
	\end{eqnarray}
	By  squaring  both sides of  \eqref{Eqq3.4} and using the  inequality ${(a+b)^2 \leq 2 (a^2+b^2)},$ for $a,b\in \mathbb R$ and then taking the expectations, one gets 
	\begin{eqnarray}
		\label{Eqq3.14}
		\mathbb E\left(\big\|\big(\widetilde f_{N,\lambda}(X_i)-f^*(X_i)\big)_i \big\|^2_n\right)&\leq& 2\, \mathbb E\left(\big\|A_N(K_N^2+\lambda I_N)^{-1}A_N^* F_N^T- F_N^T\|_n^2\right) \nonumber\\
		&&+2\, \mathbb E\left(\big\|A_N(K_N^2+\lambda I_N)^{-1}A_N^* \pmb \varepsilon^T\|_n^2\right)+4\, \mathbb E\left( \big\|\big(f^*(X_i)-\pi_N f^*(X_i)\big)_i\big\|^2_n\right)\nonumber\\
		&& +4\, \mathbb E \left( \big\|\big(\pi_N f^*(X_i)-\widetilde\pi_N f^*(X_i)\big)_i\big\|^2_n\right).
	\end{eqnarray}
	Finally, by combining the previous inequality with \eqref{Eqq3.5}, \eqref{Eqq3.6}, \eqref{SVD},  \eqref{Eqq3.11} and \eqref{Eqq3.12}, one gets the desired result \eqref{Riskerror}.\\
	
	\begin{remark}
		Although, the risk error \eqref{Riskerror} of our TKRR estimator  has some similarities with the risk error   of the RFF estimator, given for example in \cite{Avron2},  there is a fundamental difference between these two risks and their proofs. More precisely, the risk error \eqref{Riskerror} involves   the eigenvalues  of the integral operator $T_{\mathbb K},$ as well as a weighted tail trace of this latter. Moreover, it works for a fairly general random kernel matrix with unknown random sampling set. The empirical risk error of the RFF estimator involves  the eigenvalues of the kernel matrix ${\bf K}$. A similar remark is also valid for the empirical risk error of the spectrally truncated KRR estimator developed by  \cite{Amini}.
	\end{remark}
	
\begin{remark} The condition $f^*\in \mathcal H$ can be relaxed. We may assume that the RKHS $\mathcal H$ is a closed subspace of $L^2(\pmb X, d P)$ and that the regression function $f^*$ belongs to another subspace $\mathcal F$  of $L^2(\pmb X, dP).$ Let $\pi_{\mathcal H}$ be the orthogonal projection operator over $\mathcal H $ and let $g^*= \pi_{\mathcal H} f^* \in \mathcal H.$ Then by writing $f^*= g^*+(f^*-g^*)$ and using the same proof of the previous theorem, one gets
the following empirical risk error for $f^*\in \mathcal F,$
\begin{eqnarray*}\label{RiskError}
			\mathbb{E}\Big( \big\|\big(\widehat{f}_{N,\lambda}-f^*\big)(X_i)\big\|_n^2\Big)&\leq&\frac{\lambda}{2}+2\frac{\sigma^2}{n}\sum_{j=1}^N\mathbb E_X\Big[\Big(\frac{\mu_j^2}{\mu_j^2+\lambda}\Big)^2\Big]\nonumber\\
		&+& C_1\, \Big(\lambda_{N+1} \|\pi_H f^*\|^2_{\mathcal H}+\frac{1}{N}\sum_{k\geq N+1} \lambda_k \|\varphi_k\|^2_\infty  +\| f^*-\pi_{\mathcal H} f^*\|^2_{L^2(dP)}\Big).
		\end{eqnarray*}
The previous  general setting empirical risk error depends on the extra projection  error $\| f^*-\pi_{\mathcal H} f^*\|_{L^2(dP)},$ that has to be estimated. 
\end{remark}
	
	Next, we use the risk error \eqref{Riskerror} to study the two main issues related to our proposed TKRR: The choice of the convenient value of the regularization parameter $\lambda,$ and    the convergence rate of our TKRR estimator.
	Note that from  \eqref{Riskerror}, the previous two issues are essentially based on the behavior of the sequence of the eigenvalues
	$(\lambda_{k})_{k\geq 1}.$ In particular, the decay rate of the 
	$\lambda_{k}$ and the smallest rank $N_0\in \mathbb N$
	(depending only on the kernel $\mathbb K$), from which this decay rate holds true, play an important role in deriving the optimal convergence rate of our TKRR estimator. For this purpose, we consider the two frequently encountered types of the spectra:
	The spectra with a polynomial decay rate and those with an  exponential decay rate. More precisely, we assume that for some real $a\geq 0,$ the eigenfunctions $\varphi_k$ of the integral operator $T_{\mathbb K}$ satisfy the following uniform bound
	\begin{equation}\label{Eqq3.15}
		\|\varphi_k \|_\infty \lesssim k^{a/2},\qquad k\geq 1.
	\end{equation}
	Note that for some 
	positive real numbers $s > a \geq 0$ and $b>0,$ the exponential and polynomial decay rates of the  eigenvalues of $T_{\mathbb K},$  are given by 
	\begin{equation}\label{Eqq3.16}
		\lambda_{k} \lesssim e^{-b k},\quad\forall\, k\geq N_b
	\end{equation}
	and
	\begin{equation}\label{Eqqq3.16}
		\lambda_{k} \lesssim k^{-2s},\quad\forall\, k\geq N_s,
	\end{equation}
	respectively. Here, $N_b$ and $N_s$ are the ranks from which these decay rates hold true. In particular, since in \eqref{Riskerror}, we have ${\displaystyle \frac{\sigma^2}{n}\sum_{j=1}^N\mathbb E_X\left[\Big(\frac{\mu_j^2}{\mu_j^2+\lambda}\Big)^2\right]\leq \sigma^2\frac{N}{n},}$ then by using \eqref{Eqq3.15} and \eqref{Eqq3.16} for the exponential decay rate case, the inequality 
	\eqref{Riskerror} is rewritten as follows
	\begin{equation}\label{Riskerror2}
		\mathbb{E}\Big( \big\|\big(\widehat{f}_{N,\lambda}-f^*\big)(X_i)\big\|_n^2\Big)\leq \frac{\lambda}{2}+2 \sigma^2\frac{N}{n} + C_2\, N^\eta \lambda_N,\quad \eta=\max(0,a-1).
	\end{equation}
	Consequently,  the optimal values of the parameters $\lambda, N$ in terms of the  dataset  size $n,$ are given by 
	$$
	\lambda = O\Big(\sigma^2 \frac{N}{n}\Big) = O\Big( N^{\eta} e^{-b N}\Big),\quad N\geq N_b.
	$$
	Straightforward computations show that  this is the case when 
	\begin{equation}\label{Eqq3.17}
		N=\max\Big(N_b, \frac{1}{b}\log\left(\frac{n}{\sigma^2}\right)\Big),\quad \lambda=
		O\left(\sigma^2 \min\Big(\frac{N_b}{n},\frac{1}{b\, n}\log\big(\frac{n}{\sigma^2}\big)\Big)\right).  
	\end{equation}
	In a similar manner, when the spectrum of $T_{\mathbb K}$ has a polynomial decay rate and by using  \eqref{Eqqq3.16}, the optimal values of the parameters $N, \lambda$ are given by 
	\begin{equation}\label{Eqq3.18}
		N=\max\left(N_s, \Big(\frac{n}{\sigma^2}\Big)^{1/(2s+\gamma)}\Big)\right),\quad \lambda=
		O\left( \min\Big(\frac{N_s}{n},\Big(\frac{\sigma^2}{n}\Big)^{1-1/(2s+
			\gamma)}\right),\quad \gamma=\min(1,2-a).  
	\end{equation}
	We have just proved the following proposition that provides us with the optimal convergence rate of our TKRR estimator when $f^*\in \mathcal H$  and  in the presence of an exponential or a polynomial decay rate of the eigenvalues $\lambda_k.$
	
	\begin{proposition}
		Under the hypotheses of Theorem~2, the optimal convergence rates of the TKRR estimator $\widehat{f}_{N,\lambda}$ under the assumption 
		\eqref{Eqq3.16} is  given by 
		\begin{equation}\label{Riskerror2-2}
			\mathbb{E}\Big( \big\|\big(\widehat{f}_{N,\lambda}-f^*\big)(X_i)\big\|_n^2\Big)=O\left(\sigma^2 \max\Big(\frac{N_b}{n},\frac{1}{b\, n}\log\big(\frac{n}{\sigma^2}\big)\Big)\right), \quad N\geq N_b.
		\end{equation}
		Moreover, under the assumption  \eqref{Eqqq3.16}, we have 
		\begin{equation}\label{Riskerror2-3}
			\mathbb{E}\Big( \big\|\big(\widehat{f}_{N,\lambda}-f^*\big)(X_i)\big\|_n^2\Big)=O\left( \max\Big(\sigma^2\frac{N_s}{n},\Big(\frac{n}{\sigma^2}\Big)^{-\frac{(2s+\gamma-1)}{(2s+\gamma)}}\Big)\right), \quad \gamma=\min(1,2-a),\quad N\geq N_s.
		\end{equation}
	\end{proposition}

	\begin{remark}
		By using the concept of the uniform degrees of freedom at level $\epsilon >0$  of a positive and compact integral operator $T_{\mathbb K},$ one may generalize and improve  the result \eqref{Riskerror2-2} of the previous proposition.  Then, this uniform  degrees of freedom is denoted by $d_\infty(T_{\mathbb K},\epsilon)$ and defined by 
		\begin{equation}\label{degreesFreedom}
			d_\infty(T_{\mathbb K},\epsilon)= \min\{ k\in \mathbb N;\, \lambda_k \leq \epsilon\}= N_{T_{\mathbb K}}(\epsilon)=\#\{ \lambda_j;\, \lambda_j >\epsilon \}+1.
		\end{equation}
		In this case,  a  refined truncation order $N,$ appearing in the right hand-side of \eqref{Riskerror2} is simply given by 
		\begin{equation}\label{refinedN}
			N=N(\epsilon_n)=d_\infty(T_{\mathbb K},\epsilon_n)\quad\mbox{with}\quad\epsilon_n =\frac{\sigma^2}{n}\big(d_\infty(T_{\mathbb K},\epsilon_n)\big)^{1-\eta},\quad \eta=\max(0,a-1). 
		\end{equation} 
	\end{remark}
	
	Besides its good convergence  properties, given by the previous proposition, our TKRR estimator has the desirable property of requiring a low  computational load, compared with the full KRR estimator. 
	
	\begin{remark}
		From \eqref{GramMatrix1} and \eqref{Estimator1}, the TKRR estimator has time complexity of  $O(n N^2),$ whereas the  full KRR estimator requires time complexity of  $O(n^3).$ Since, in general $N\ll n,$ then the TKRR is much faster than the classical KRR. On the other hand, our TKRR has the same time complexity as the spectrally truncated KRR, given in \cite{Amini}. Nonetheless our TKRR does not require any SVD decomposition and handles random sampling dataset drawn from a fairly general unknown multivariate probability law. Moreover, the  proposed TKRR is competitive with   the combined sketching-KRR estimator, given in \cite{Yang}. This last estimator requires a time complexity of $O(n^2 N^2)$ for dense sketches and an $O(n^2 \log(N)),$ for Fourier and Hadamard sketches.   
	\end{remark}
	
	An other important feature of our TKRR is its optimal  convergence rate that  coincides with the optimal convergence rate of the full KRR estimator.  More precisely, the following proposition shows that  the full KRR estimator has the same convergence rates \eqref{Riskerror2-2} and \eqref{Riskerror2-3}. Note that the full KRR estimator $\widehat f_n^\lambda$ is given by \eqref{Estimator3}--\eqref{Tikhonov4}.  
	
	\begin{proposition}
		Under the hypotheses of Theorem~1 and Theorem~2 and  the  assumptions \eqref{Eqq3.16} and \eqref{Eqqq3.16}, the optimal convergence rate of the full KRR estimator $\widehat f_n^\lambda$  coincide with the convergence rates \eqref{Riskerror2-2} and \eqref{Riskerror2-3}, respectively.
	\end{proposition}
	
	\noindent
	{\bf Proof:} To derive the empirical risk error of the full KRR estimator \eqref{Estimator3}--\eqref{Tikhonov4}, one might use a  simple adaptation of the proof of Theorem 2. More precisely, instead of considering the pseudo-inverse of the regularized truncated matrix $A_N,$ given by \eqref{Eqq3.3},  it suffices to consider the inverse of the regularized full Gram matrix 
	${\displaystyle A_n +\lambda I_n},$ given by \eqref{Tikhonov3}.  Also, in this case and thanks to the representer theorem and the hypothesis that $f^*\in \mathcal H$ and as it is done in \cite{Amini}, we may assume that $f^*(X_i)= \widetilde \pi_n f^*(X_i)$ for $1\leq i\leq n.$
	In this case, the empirical risk error of $\widehat f_n^\lambda$ is given by 
	\begin{equation}\label{Riskerror4}
		\mathbb{E}\Big( \big\|\big(\widehat f_n^\lambda-f^*\big)(X_i)\big\|_n^2\Big)\leq \frac{\lambda}{2}+2\frac{\sigma^2}{n}\sum_{j=1}^n\mathbb E_X\Big[\Big(\frac{\lambda_j(A_n)}{\lambda_j(A_n)+\lambda}\Big)^2\Big]
	\end{equation}
	Since ${\displaystyle 0<\frac{\lambda_j(A_n)}{\lambda_j(A_n)+\lambda} <1,}$ then for any integer $1\leq N<n,$ we have 
	\begin{eqnarray}\label{Eqq3.19}
		\mathbb{E}\Big( \big\|\big(\widehat f_n^\lambda-f^*\big)(X_i)\big\|_n^2\Big)&\leq& \frac{\lambda}{2}+2 \sigma^2\frac{N}{n}+2\frac{\sigma^2}{n}\sum_{j=N+1}^n\mathbb E_X\Big(\frac{\lambda_j(A_n)}{\lambda_j(A_n)+\lambda}\Big)\nonumber\\
		&\leq & \frac{\lambda}{2}+2\sigma^2\frac{N}{n}+2\frac{\sigma^2}{n\,\lambda}\sum_{j=N+1}^n\mathbb E_X\big(\lambda_j(A_n)\big).
	\end{eqnarray}
	Under the hypothesis that the sequence of the eigenvalues of the integral operator $T_{\mathbb K}$ has an exponential decay, given by \eqref{Eqq3.16} and by using the previous inequality and estimate \eqref{compare2},  one gets for $N_b\leq N <n,$
	\begin{equation}\label{Eqq3.20}
		\mathbb{E}\Big( \big\|\big(\widehat f_n^\lambda-f^*\big)(X_i)\big\|_n^2\Big)\lesssim \frac{\lambda}{2}+2 \sigma^2\frac{N}{n}+2\frac{\sigma^2}{n\,\lambda} \lambda_{N+1} \|\varphi_{N+1}\|^2_\infty.
	\end{equation}
	Consequently, the optimal values of $\lambda$ and $N$ satisfy the  estimates
	\begin{equation}\label{Eqq3.21}
		\lambda = O\Big(\sigma^2\frac{N}{n}\Big)= O\Big(\frac{\sigma^2}{n} N^\eta e^{-bN}\Big),\quad \eta=\max(0,a-1),\quad N\geq N_b.
	\end{equation}
	By combining \eqref{Eqq3.20} and \eqref{Eqq3.21}, one concludes that under hypothesis \eqref{Eqq3.16}, our TKRR estimator has the same convergence rate as the full KRR.
	In a similar manner, one can easily check that this is also the case under the hypothesis \eqref{Eqqq3.16}. We leave the details for the reader. 
	
	\section{Examples and numerical simulations}
	
	In  this  section, we first  illustrate the results of section 3 and section 4. We consider two   examples of the Sinc and   Gaussian kernels. These  Mercer's kernels are frequently used in the framework of KRR based NP regression estimators. Then, we give numerical simulations that illustrate the theoretical properties  of our proposed TKRR estimator.\\
	
	\subsection {The Sinc and Gaussian kernels examples}
	\subsubsection{Sinc kernel case}
	It is well known that there is a rich literature concerning the spectral analysis of the Sinc kernel operator, as well as the various
	mathematical statistics and signal processing applications related to this kernel, see for example \cite{Hogan}. This kernel and its associated integral operator are briefly described as follows. Consider a positive real number $c>0$ called bandwidth, then the uni-dimensional Sinc-kernel is defined by  
	\begin{equation}
		\label{SincKernel}
		\mathbb{K}_c(x,y)=\frac{\sin(c(x-y))}{\pi(x-y)},\quad x,y\in I=[-1,1].
	\end{equation}
	Since $\mathbb{K}_c(x,y)=\kappa(x-y),$ where $\kappa(\cdot)$ is the Fourier transform of the uniform probability measure over $I=[-1,1],$ then by Bochner's theorem, $\mathbb{K}_c(\cdot,\cdot)$ is a positive definite kernel. It is well known that  the associated  RKHS is given by the space of  bandlimited functions with bandwidth $c.$ That is the subspace of functions from $L^2(\mathbb R)$ with Fourier transforms supported in $[-c,c].$ The usual probability measure associated with the Sinc-kernel is given by the uniform measure ${\displaystyle d \mu(x)= \frac{1}{2} \mathbf 1_{[-1,1]}(x) dx.}$
	We recall that if $\mathcal Q_c$ is the Hilbert-Schmidt operator defined on $L^2(I,d\mu)$ with kernel $\mathbb{K}_c(x,y)$, then we have  $$\mathcal{Q}_c(\varphi_{m,c})(x)=\int_{-1}^{1} \frac{\sin(c(x-y))}{\pi(x-y)}\varphi_{m,c}(y)dy=\lambda_m(c) \varphi_{m,c}(x), \quad x\in [-1,1].$$
	Here, $\varphi_{m,c},\,  \lambda_m(c)$ is the $m-$th eigenfunction and the associated eigenvalue of $\mathcal{Q}_c.$ The $\varphi_{m,c}$ are known as the prolate spheroidal wave functions (PSWFs). The rich properties of the PSWFs and the behavior of  their eigenvalues have been extensively studied in the literature since the pioneer work in the subject starting from  the 1960's by  D. Slepian and his co-authors H. Landau and H. Pollak, see \cite{Hogan} for more details.  Since these early works on the PSWFs, it is known that the sequence of the $\lambda_m(c)$ has an asymptotic super-exponential decay rate. Recently, it has been shown in  \cite{Bonami-Karoui}, that the  optimal asymptotic  super-exponential decay rate
	of the $\lambda_m(c)$ is given by 
	$$\lambda_m(c) \lesssim \exp\Big(-2m \log\Big(\frac{4 m}{e c}\Big)\Big),\quad m\gg 1.$$
	Also, in  \cite{Bonami-Karoui2},  the following non-asymptotic decay rate of the $\lambda_m(c)$ has been given,
	\begin{equation}
		\lambda_m(c) \leq \exp\Big(-(2m+1) \log\Big(\frac{2(m+1)}{e c}\Big)\Big),\quad m\geq \frac{ec}{2} =N_c.
	\end{equation}
	Consequently, for $d=1$ and for  the Sinc kernel case, one can take $b=2$ and ${\displaystyle N_b=\frac{ec}{2}}$ in the empirical risk error \eqref{Riskerror2-2} of Proposition 1. That is 
	\begin{equation}\label{Riskerror2-4}
		\mathbb{E}\Big( \big\|\big(\widehat{f}_{N,\lambda}-f^*\big)(X_i)\big\|_n^2\Big)=O\left(\sigma^2 \max\Big(\frac{e c}{2 n},\frac{1}{2\, n}\log\big(\frac{n}{\sigma^2}\big)\Big)\right), \quad N\geq \frac{ec}{2}.
	\end{equation}
	Note that for large values of the bandwidth $c,$ the decay rank $N_c=\frac{e c}{2}$ is pessimistic. In fact, it is  known from the work of H. Landau,  that  the decay region of the $\lambda_n(c)$ starts at ${\displaystyle N_{1,c}= \frac{2c}{\pi} + O\big(\log(c)\big).}$ More precisely, it has been shown that for the Sinc kernel, the uniform degrees of freedom at level $\epsilon >0$ is given by 
	$$ d_{\infty}(T_{\mathbb K_c},\epsilon)= \frac{2c}{\pi} +\frac{1}{\pi^2} \log\Big(\frac{1-\epsilon}{\epsilon}\Big)\log\Big(\frac{2c}{\pi}\Big)+ o\big(\log c\big),\quad c\gg 1.$$
	Also, it has been shown in \cite{Bonami-Karoui2} that 
	\begin{equation}\label{2.2.6}
		\sup_{x \in [-1,1]}|\varphi_{m,c}(x)|^2\lesssim \big(m(m+1)+c^2\big)^\frac{1}{2} \lesssim m, \quad \forall \; m\geq  \Big[ \frac{2c}{\pi}\Big].
	\end{equation}
	That is by using \eqref{refinedN} for the Sinc kernel with $a=1$ and the previous uniform degrees of freedom, a refined truncation order is given by 
	\begin{equation}\label{refinedorder1}
		N=N(\epsilon_n)=\left\lceil\frac{2c}{\pi} +\Big(\frac{1}{\pi^2} \log\big(\frac{1}{\epsilon_n}\big)+1\Big)\log\Big(\frac{2c}{\pi}\Big)
		\right\rceil,
	\end{equation}
	where $\epsilon_n$ is a solution of the equation
	\begin{equation}\label{refinedorder2}
		\epsilon_n = \frac{\sigma^2}{n} \left(\frac{2c}{\pi} +\Big(\frac{1}{\pi^2} \log\Big(\frac{1}{\epsilon_n}\Big)+1\Big)\log\Big(\frac{2c}{\pi}\Big)
		\right).
	\end{equation}
	Here, $\lceil x \rceil$ denotes the smallest integer greater or equal  to $x.$ Finally, for the dimension $d\geq 2,$ one may use the tensor product Sinc kernel 
	$$\mathbb K_c^d(\pmb x,\pmb y)=\prod_{i=1}^d \mathbb K_c(x_i-y_i),\quad \pmb x= \big(x_i\big)_{i=1}^d,\quad \pmb y= \big(y_i\big)_{i=1}^d\in I^d.$$
	The  $d-$variate eigenfunctions associated to  $\mathbb K_c^d(\pmb x,\pmb y)$ and their corresponding eigenvalues are given by the tensor products of the $\varphi_{k,c}(x_i)$ and their corresponding eigenvalues $\lambda_k(c).$ This allows us to derive  a decay rate estimate for the eigenvalues, as well as a uniform degrees of freedom associated with a $d-$dimensional  tensor product of Sinc kernels. We leave the details for the readers. 
	
	\subsubsection{The Gaussian  kernel case} As for the sinc kernel, we first recall  the uni-dimensional Gaussian kernel ${\displaystyle 
		\mathbb K_G(x,y),\, x,y\in\mathbb R}$ and some of its associated spectral properties. Then, we show how to extend these 
	properties to the $d-$dimensional case. It is well known that the Gaussian kernel is one of the most used kernels in mathematical statistics and machine learning applications. This kernel is described as follows. Consider two positive real  numbers $\xi, c >0,$ then  $\mathbb K_G$ is the positive definite  convolution  kernel with associated usual probability measure $d P_c(x),$  given by  
	\begin{equation}
		\label{GaussianKernel}
		\mathbb K_G(x,y)= e^{- \xi (x-y)^2},\qquad dP_c(x)=\sqrt{\frac{\pi}{2c}} e^{- 2 c x^2}\quad x,y\in \mathbb R.
	\end{equation}
	From \cite{Zhu2}, the eigenvalues of the integral operator $T_{\mathbb K_G}$ are given by
	\begin{equation}\label{eigenvalsG}
		\lambda_k = \sqrt{\frac{\pi}{\xi+c+\sqrt{c^2+2c\xi}}} \left(\frac{\xi}{c+\xi+\sqrt{c^2+2c\xi}}\right)^{k-1},\quad k\geq 1.
	\end{equation}
	The associated eigenfunctions  $\phi_k,\, k\geq 1$ are the well known dilated Hermite functions, given by 
	\begin{equation}\label{eigenfctsG}
		\varphi_k(x,\gamma) = \frac{1}{\sqrt{2\gamma}} \phi_k(\sqrt{2\gamma} x)= \frac{1}{\sqrt{2\gamma}} \alpha_{k-1} H_{k-1}\big(\sqrt{2 \gamma} x\big)e^{- \gamma x^2},\quad \gamma=\sqrt{c^2+2c\xi},\quad k\geq 1.
	\end{equation}
	Here, ${\displaystyle \alpha_k= \frac{2^{-k/2}}{\pi^{1/4} \sqrt{k!}}}$ and ${\displaystyle H_k(x)= (-1)^k e^{x^2} \frac{d^k}{d x^k}\Big(e^{-x^2}\Big),\, k\geq 0}$ are the Hermite polynomials. 
	Moreover, it is well known that the normalized Hermite functions $\phi_k$ satisfy the following uniform bound,
	\begin{equation}\label{eigenfcts2G}
		\sup_{x \in \mathbb R} |\phi_k(x)|\leq 1,\quad k\geq 1.
	\end{equation}
	Hence, by using \eqref{eigenvalsG}--\eqref{eigenfcts2G}, as well as the empirical risk error \eqref{Riskerror2} with $\eta=0,$ one can get the optimal values of the truncation order $N$ and the regularization parameter $\lambda>0,$ in terms of the parameters $\xi, \, c >0$ and according to the rule \eqref{Eqq3.17} with ${\displaystyle b=\frac{c+\xi+\sqrt{c^2+2c\xi}}{\xi}.}$

	 Finally, for the general $d-$dimensional case, the Gaussian kernel and its associated probability measure are  given by 
	$$\mathbb K_G^d(\pmb x,\pmb y)= e^{- \xi \|\pmb x-\pmb y\|^2_{2}},\quad dP_c(\pmb x)=\left(\frac{\pi}{2c}\right)^{d/2} e^{- 2 c \|\pmb x^2\|_2^2},\quad  \pmb x, \pmb y \in \mathbb R^d.$$
	The associated $d-$variates eigenfunctions and their corresponding eigenvalues are given by the tensor products of the uni-variate
	$\varphi_k(x_i,\gamma)$ and their  corresponding $\lambda_k.$
	
	\subsection{Numerical simulations}
	
	In this paragraph, we give some numerical simulations that illustrate the different results of this work, in particular the results of  Theorem 1 and Theorem 2 and the  consequences of this latter. We restrict ourselves to the uni-dimensional case, since the aim is to illustrate the theoretical results of this work.  These numerical simulations are given by the following three examples.\\
	
	\noindent{\bf Example 1:} In this first example, we illustrate the results of Theorem~1 and its corollary 1. For this purpose, we have first considered the Sinc-kernel $\mathbb K_c,$ given by \eqref{SincKernel} with $c=25.$  Then, instead of using its associated usual sampling probability measure $dP(x)=\frac{1}{2} \mathbf 1_{[-1,1]}(x) dx,$ we have considered a set of $n=200$ i.i.d samples following a truncated standard normal distribution, supported on $I=[-1,1]$ with probability density function 
	\begin{equation}
		\label{trincatedNormal}
		\rho(x)= \frac{\phi(x)}{\Phi(1)-\Phi(-1)} {\mathbf 1}_{[-1,1]}(x),\quad \phi(x)=\frac{1}{\sqrt{2\pi}} \exp\Big(-\frac{x^2}{2}\Big),\quad \Phi(x)=\frac{1}{2}\big(1+\mbox{erf}(x/\sqrt{2})\big).
	\end{equation}   
	It is well known that the integral operator $T_{\mathbb K_c}=\mathcal Q_c$ is one to one. Moreover, it is easy to see  that $d\rho(x)$ is absolutely continuous with respect to the uniform probability measure $dP(x).$
	Hence, the conditions of Theorem 1 are satisfied for these two probability measures. We have  considered the truncation order $N=25$ and computed the average over $10$ realizations for the eigenvalues of the full Gram matrix $B^c_n$ and its main $n\times N$ sub-matrix $A_N^c,$ given by 
	$${\displaystyle  B^c_n= \frac{1}{n}\left[\mathbb K_c(X_i, X_j)\right]_{1\leq i,j\leq n},\qquad\,\,  A^c_N= \frac{1}{n}\left[\mathbb K_c(X_i, X_j)\right]_{1\leq i\leq n, 1\leq j\leq N}}.$$
	Recall that for each realization, these Gram matrices are constructed by applying the Sinc-kernel on $n=200$ i.i.d random sampling points drawn from $\rho(\cdot).$ In Figure 1(a), we have plotted the graphs of the  averages $\overline{\lambda_j}(B^c_n)$ and $\overline{\mu_j}(A^c_N)$ of   the eigenvalues $\lambda_j(B^c_n),$ and the singular values $\mu_j((A^c_N)),$ versus the true eigenvalues $\lambda_j(T_{\mathbb K_c})$ of the Sinc-kernel integral operator $T_{\mathbb K_c}=\mathcal Q_c.$ Also, to capture the fast decay rate of the previous sequences of eigenvalues, we have plotted in Figure 1(b),  the graphs of the logarithms of the averages of these  eigenvalues. Moreover, to illustrate the decay estimates for the trace tails $\mathcal T_k(B^c_n)$ and  $\mathcal T_k(A^c_N)$ of $B^c_n$  and $A^c_N,$ given by \eqref{compare2} and \eqref{SingularValues1}, we have plotted in Figure 2(a) the graphs of $\mathcal T_k(B^c_n)$ and $\mathcal T_k(A^c_N)$ versus the true eigenvalues $\lambda_k(T_{\mathbb K_c}),$ $1\leq k\leq N.$ In figure 2(b), we have plotted the graphs of  $\log(\mathcal T_k(B^c_n)), \log(\mathcal T_k(A^c_N))$ versus the true eigenvalues $\log(\lambda_k(T_{\mathbb K_c})),$ for different values of $1\leq k\leq N.$ The numerical simulations given by Figure 1 and Figure 2 are   coherent with the theoretical results of Theorem 1 and corollary 1. \\
	
	Also,  we have considered the Gaussian kernel $K_G(x,y)=e^{- \xi (x-y)^2},\, x,y\in I=[-1,1]$ with $\xi=25.$   The associated 
	Gram matrix $B^G_n$ and its main $n\times N$ sub-matrix $A^G_N,$ for $ N<n$ are given by
	$${\displaystyle  B^G_n= \frac{1}{n}\left[\mathbb K_G(X_i, X_j)\right]_{1\leq i,j\leq n},\qquad\,\,  A^G_N= \frac{1}{n}\left[\mathbb K_G(X_i, X_j)\right]_{1\leq i\leq n, 1\leq j\leq N}}.$$
	Here, the $X_i$ are i.i.d random samples drawn from the uniform law over $I=[-1,1].$ Although  
	the  uniform probability measure is not the usual Gaussian measure for the kernel $\mathbb K_G$  and as predicted by  Theorem~1 and its Corollary~1, the  eigenvalues of $B^G_n$ and $A_N^G$  have similar decay rate as the true eigenvalues sequence $\lambda_k,$ given by \eqref{eigenvalsG}, with $\xi=25$ and $c=1.$ In Figure 3, we plotted the  eigenvalues of the integral operator $\lambda_j(T_{\mathbb K_G}),$ according to \eqref{eigenvalsG}, as well as $\overline{\lambda_j}(B^G_n),\, \overline{\mu_j}(A^G_N),$ the averages over $10$ realizations for the eigenvalues of $B^G_n$ and the singular values of $A_N^G,$ for $n=200$ and $N=25.$

	\begin{figure}[]\hspace*{0.0cm}
		{\includegraphics[width=15.05cm,height=5.5cm]{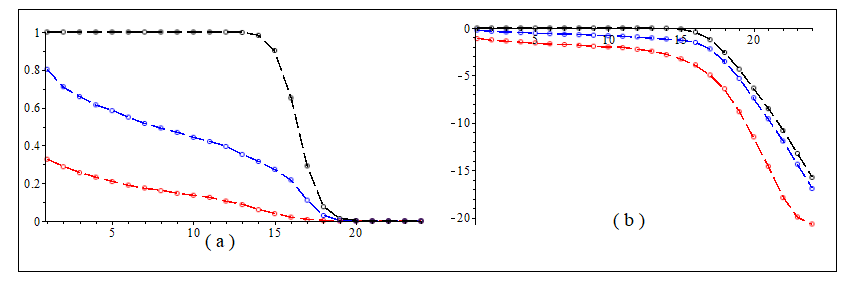}}
		\vskip -0.5cm\hspace*{1cm} \caption{(a) Graphs of  $ \lambda_j(T_{\mathbb K_c})$ (black) versus $\overline{\lambda_j}(B^c_n) $ (blue) and the singular values $\overline{\mu_j}(A^c_N)$ with $c=25,$ $n=200,$ $N=25$.
			(b) Graphs of   $\log(\lambda_j(T_{\mathbb K_c}))$ (black) versus $\log(\overline{\lambda_j}(B^c_n)) $ (blue) and  $\log(\overline{\mu_j}(A^c_N))$ (red). }
	\end{figure}
	\begin{figure}[]\hspace*{0.0cm}
		{\includegraphics[width=15.05cm,height=5.5cm]{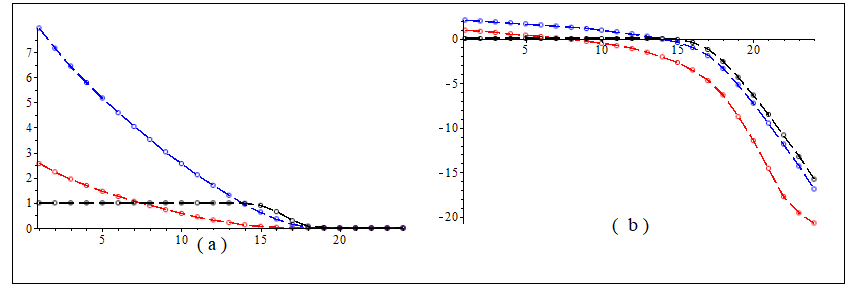}}
		\vskip -0.5cm\hspace*{1cm} \caption{(a) Graphs of  the trace  tails $\mathcal T_k(B_n^c)$ (blue) and $\mathcal T_k(A_N^c)$ (red), versus $\lambda_k(T_{\mathbb K_c})$ (black), c=$25,$ $n=200,$ $N=25.$  
			(b) Graphs of   $\log(\mathcal T_k(B_n^c))$ (blue) and  $\log(\mathcal T_k(A_N^c)) $ (red), versus    $\log(\lambda_k(T_{\mathbb K_c}))$ (black). }
	\end{figure}
	\begin{figure}[]\hspace*{0.0cm}
		{\includegraphics[width=15.05cm,height=5.5cm]{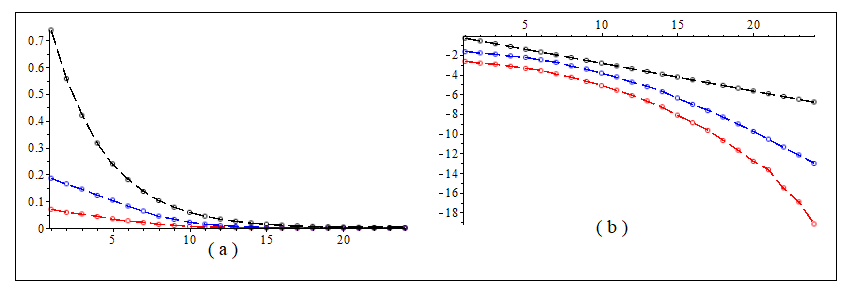}}
		\vskip -0.5cm\hspace*{1cm} \caption{(a) Graphs of  $ \lambda_j(T_{\mathbb K_G})$ (black) versus $\overline{\lambda_j}(B^G_n) $ (blue) and the singular values $\overline{\mu_j}(A_N^G)$ with $n=200,$ $N=25$, $\xi=25, c=1.$
			(b) Graphs of   $\log(\lambda_j(T_{\mathbb K_G}))$ (black) versus $\log(\overline{\lambda_j}(B^G_n)) $ (blue) and  $\log(\overline{\mu_j}((A_N^G)))$ (red). }
	\end{figure}
	
	\vskip 1cm
	\noindent{\bf Example 2:}  In this second example, we illustrate the results of Theorem 2, as well as Proposition 1 and Proposition 2. For this purpose, we  consider the non-parametric regression problem \eqref{model1} with $d=1$ and the  true regression function $f^*,$
	given by $${\displaystyle f^*(x)=\frac{\sin (20\, x)}{20\, x},\,\, x\in I=[-1,1]}.$$ We use the Mercer's kernel $\mathbb K=\mathbb K_c,$ the Sinc-kernel with the special value of $c=25$. It is well known that when associated with the uniform measure over $I,$  the  RKHS space $\mathcal H,$  associated with $\mathcal K_c$ is the space of functions from  $L^2(\mathbb R)$ with Fourier transforms supported in $[-c,c].$ Hence, the previous test function $f^* \in \mathcal H.$ We have considered $n=200$ i.i.d $X_i$ drawn from the truncated normal distribution of the previous example. Moreover, the added i.i.d noises $\varepsilon_i$ are drawn from a centered normal distribution with the two different values of standard deviations $\sigma=0.1$ and $\sigma=0.5.$ Then, we have constructed our estimator \eqref{Estimator1} with   the regularization parameter $\lambda,$ chosen according to the rule \eqref{Eqq3.17} with 
	$b=2$ and $N_b=\Big\lceil\frac{ec}{2}\Big\rceil=34,$ that is $\lambda=2.47e-4.$  We have considered the four values of truncation order $N=20, 25, 30, 50.$ We  computed 
	the mean over $10$ realizations for the empirical risk errors $\mathcal{\widehat R}_N$ and the theoretical empirical risk error $\mathcal R_n,$ given by \eqref{Riskerror2-4}. That is 
	$$\mathcal{\widehat R}_N = \frac{1}{n} \sum_{i=1}^n \Big( \widehat f_{N,\lambda} (X_i)-f^*(X_i)\Big)^2,\quad \mathcal R_n = O\Big(\sigma^2 \frac{e c}{2\, n}\Big).$$
	The obtained numerical results are given by Table 1.  These results indicate that the theoretical empirical risk error bounds given by Theorem 2 and Proposition 1 are fairly tight. Moreover and as we have previously mentioned, the refined optimal truncation order for the Sinc-kernel case is given by \eqref{refinedorder1}-\eqref{refinedorder2}. In particuler, for $c=25,$ this optimal truncation order is given by $N_1=21$ for $\sigma=0.1$ and $N_2=20$ for $\sigma=0.5.$ This explains why the empirical risk errors of Table 1, that are  obtained for $N=20$ are similar to those obtained for larger values of $N.$ \\

	\begin{center}
		\begin{table}[h]
			\vskip 0.2cm\hspace*{2.0cm}
			\begin{tabular}{cccccccc} \hline
				& & & & & & & \\
				$\sigma$   &$N$&$\mathcal{\widehat R}_N$&$\mathcal R_n$        &   $\sigma$  &$N$&$\mathcal{\widehat R}_N$&$\mathcal R_n$   \\   \hline
				$0.1$& $20$ & $9.05e-3  $      &$1.70e-3  $ & $0.5$ & $ 20$& $2.07e-2  $& $4.25e-2     $ \\
				$-$& $25$ & $4.64e-3  $      &$1.70e-3  $ & $-$ & $ 25$& $2.21e-2  $& $4.25e-2     $ \\
				$-$& $30$ & $1.02e-3  $      &$1.70e-3  $ & $-$ & $ 30$& $2.20e-2  $& $4.25e-2     $ \\
				$-$& $50$ & $1.00e-3  $      &$1.70e-3  $ & $-$ & $ 50$& $1.92e-2  $& $4.25e-2     $ \\

				\\ \hline
			\end{tabular}
			\caption{Numerical results of example 2.}
			\label{tableau2}
		\end{table}
	\end{center} 
	
	\vskip 1cm
	\noindent{\bf Example 3:}  In this last example, we consider the Gaussian kernel $K_G(x,y)=e^{- \xi (x-y)^2},\, x,y\in I=[-1,1]$ with $\xi=25$  and $c=1.$   Then, we consider similar  numerical simulations as the previous example, but with the  new  synthetic regression function
	\begin{equation}
		f_{\xi,c}^*(x)= \exp\Big(-\sqrt{c^2+c\xi}\Big) \cdot \Big(1+\sum_{j=1}^{10} \frac{x^j}{j}\Big),\quad \xi=25,\, c=1.
	\end{equation}
	Note that when associated  with the Gaussian probability measure $dP_c,$ given by \eqref{GaussianKernel}, the RKHS $\mathcal H$ corresponding to $\mathbb K_G$ is spanned by the Hermite functions $\phi_k,\, k\geq 1,$ given by \eqref{eigenfctsG}. Consequently,  the previous test regression function $f_{\xi,c}^* \in \mathcal H.$ Also, at each of the $10$ realizations, we consider a new set of  
	$n=200$ i.i.d random sampling points $X_i$  drawn from the  uniform law over $I.$ In Table 2, we have listed the  mean of the empirical risk errors $\mathcal{\widehat R}_N$ and  the theoretical empirical risk errors $\mathcal R_n,$     given by \eqref{Riskerror2-2} with $N_b=20$ and $b=\frac{c+\xi+\sqrt{c^2+2c\xi}}{\xi}\approx \frac{4}{3}.$ That is
	$$\mathcal{\widehat R}_N = \frac{1}{n} \sum_{i=1}^n \Big( \widehat f_{N,\lambda} (X_i)-f_{\xi,c}^*(X_i)\Big)^2,\quad \mathcal R_n = O\Big(\frac{\lambda}{2} +\frac{2\sigma^2 N_b}{n}+\lambda_N\Big).$$
	Here, $\lambda$ and $\lambda_N$ are given by  \eqref{Eqq3.17} and   \eqref{eigenvalsG}, respectively. Again, the numerical results of Table 2 are highly coherent with the theoretical results of Theorem 2 and Proposition 1. 
	\begin{center}
		\begin{table}[h]
			\vskip 0.2cm\hspace*{2.0cm}
			\begin{tabular}{cccccccc} \hline
				& & & & & & & \\
				$\sigma$   &$N$&$\mathcal{\widehat R}_N$&$\mathcal R_n$        &   $\sigma$  &$N$&$\mathcal{\widehat R}_N$&$\mathcal R_n$   \\   \hline
				$0.1$& $20$ & $1.13e-2  $      &$4.75e-3  $ & $0.5$ & $ 20$& $1.20e-2  $& $5.01e-2     $ \\
				$-$& $25$ & $5.26e-3  $      &$4.75e-3  $ & $-$ & $ 25$& $1.43e-2  $& $5.01e-2     $ \\
				$-$& $30$ & $1.91e-3  $      &$4.75e-3  $ & $-$ & $ 30$& $1.04e-2  $& $5.01e-2     $ \\
				$-$& $50$ & $1.10e-3  $      &$4.75e-3  $ & $-$ & $ 50$& $8.72e-3  $& $5.01e-2     $ \\
				& & & & & & & \\ \hline
			\end{tabular}
			\caption{Numerical results of Example 3.}
		\end{table}
	\end{center}

	\vskip 1cm
	\noindent
	{\large{\bf Acknowledgements}} The authors would like to thank Professor Aline Bonami for the fruitful discussions and suggestions concerning the results of Theorem 1 and their possible extensions.\\

\end{document}